\documentclass[12pt,a4paper,reqno]{amsart}
\usepackage[english]{babel}
\usepackage{amssymb,latexsym,amsfonts,amsthm,upref,amsmath}
\usepackage[margin=1in]{geometry}
\usepackage[dvipsnames,x11names]{xcolor}
\definecolor{myblue}{HTML}{003399}
\usepackage{hyperref}
\hypersetup{colorlinks,citecolor=myblue,filecolor=black,linkcolor=myblue,urlcolor=myblue}
\usepackage{comment} 
\usepackage{enumerate} 
\usepackage{tikz}
\usepackage{float}
\usepackage{cleveref}
\usepackage{mathtools}
\usepackage{cite}
\usepackage{filecontents}
\makeatletter
\newcommand{\leqnomode}{\tagsleft@true}
\newcommand{\reqnomode}{\tagsleft@false}
\newcommand{\cev}[1]{\reflectbox{\ensuremath{\vec{\reflectbox{\ensuremath{#1}}}}}}
\makeatother
\newtheorem*{corintro*}{Corollary}
\newtheorem*{thm*}{Theorem}
\newtheorem*{lem*}{Lemma}
\newtheoremstyle{prim}{}{}{\normalfont}{}{\bfseries}{.}{ }{}
\newtheoremstyle{stil}{}{}{\slshape}{}{\bfseries}{.}{ }{}
\theoremstyle{stil}
\newtheorem{thm}{Theorem}[section]
\newtheoremstyle{defi}{}{}{}{}{\bfseries}{.}{ }{}
\theoremstyle{defi}
\newtheorem{defn}[thm]{Definition}
\theoremstyle{defi}
\newtheorem{rem}[thm]{Remark}
\theoremstyle{stil}
\newtheorem{pro}[thm]{Proposition}
\theoremstyle{stil}
\newtheorem{lem}[thm]{Lemma}
\theoremstyle{stil}
\newtheorem{kor}[thm]{Corollary}
\theoremstyle{prim}

\newenvironment{prf}{\noindent \textit{Proof.}}{\null\hfill$\qed$\hskip
2mm\vskip 2mm}

\newcommand{\DY}{ {\rm DY}}
\newcommand{\I}{ {\rm I}}

\newcommand{\J}{ {\rm J}}

\newcommand{\A}{ {\rm A}}
\newcommand{\V}{\mathcal{V}}
\newcommand{\W}{ {\rm W}}

\newcommand{\F}{ {\rm F}}

\newcommand{\Ft}{ \widetilde{\rm F}}
\newcommand{\Y}{ {\rm Y}}

\newcommand{\R}{ {\overline{R}}}

\newcommand{\vac}{\mathop{\mathrm{\boldsymbol{1}}}}
\newcommand{\tz}{\tau}

\newcommand{\gl}{\mathfrak{gl}}
\newcommand{\sll}{\mathfrak{sl}}

\newcommand{\CC}{\mathbb{C}}

\newcommand{\ZZ}{\mathbb{Z}}
\newcommand{\TT}{\mathbb{T}}

\newcommand{\Tf}{\mathfrak{T}}

\newcommand{\Sc}{\mathcal{S}}
\newcommand{\mc}{\mathcal{M}}
\newcommand{\Tc}{\mathcal{T}}

\newcommand{\Vc}{\mathcal{V}}

\newcommand{\wtld}{\widetilde}
\newcommand{\wht}{\widehat}
\newcommand{\wvr}{\overline}
\newcommand{\wndr}{\underline}

\newcommand{\ot}{\otimes}
\newcommand{\ts}{\hspace{1pt}}

\newcommand{\tr}{ {\rm tr}}
\newcommand{\sgn}{ {\rm sgn}}
\newcommand{\xpan}{\mathop{\mathrm{span}}}
\newcommand{\ndo}{\mathop{\mathrm{End}}}

\newcommand{\rez}{\mathop{\mathrm{Res}}}
\newcommand{\Sym}{\mathfrak S}

\newcommand{\cdotrl}{\mathop{\hspace{-2pt}\underset{\text{RL}}{\cdot}\hspace{-2pt}}}
\newcommand{\cdotlr}{\mathop{\hspace{-2pt}\underset{\text{LR}}{\cdot}\hspace{-2pt}}}

\newcommand{\fand}{\quad\text{and}\quad}
\newcommand{\Fand}{\qquad\text{and}\qquad}

\newcommand{\non}{\nonumber}
\newcommand{\beq}{\begin{equation}}
\newcommand{\eeq}{\end{equation}}
\newcommand{\ben}{\begin{equation*}}
\newcommand{\een}{\end{equation*}}

\begin{document}

\title[Quantum current algebras associated with rational $R$-matrix]{Quantum current algebras associated with rational $R$-matrix}
\author{Slaven Ko\v{z}i\'{c}} 
\address{ Department of Mathematics, Faculty of Science, University of Zagreb, 10000 Zagreb, Croatia \and School of Mathematics and Statistics F07, University of Sydney, NSW 2006, Australia}
\email{kslaven@math.hr}


\begin{abstract}
We study quantum current algebra $\textrm{A}(\overline{R})$ associated with 
the
rational $R$-matrix 
of $\gl_N$
and we give explicit formulae for the elements of its center at the critical level. Due to Etingof--Kazhdan's construction, the level $c$ vacuum module $\mathcal{V}_c(\overline{R})$ for the algebra $\textrm{A}(\overline{R})$ possesses a quantum vertex algebra structure for any complex number $c$. We prove that any module for the quantum vertex algebra $\mathcal{V}_c(\R)$ is naturally equipped with a structure of restricted  $\textrm{A}(\overline{R})$-module of level $c$ and vice versa.
\end{abstract}

\maketitle

\allowdisplaybreaks

\section*{Introduction}
Let $\mathfrak{g}$ be a Lie algebra over $\CC$ equipped with a symmetric invariant bilinear form and let $\widehat{\mathfrak{g}}=\mathfrak{g}\ot\CC[t,t^{-1}] \oplus \CC C$ be the corresponding affine Lie algebra. For any complex number $c$ we associate with $\widehat{\mathfrak{g}}$ the induced module
$$V_c (\mathfrak{g}) = U(\widehat{\mathfrak{g}})\ot_{U(\widehat{\mathfrak{g}}_{(\leqslant 0 )})} \CC_c,\quad\text{where}\quad \widehat{\mathfrak{g}}_{(\leqslant 0 )}=\coprod_{n\leqslant 0}\left( \mathfrak{g}\ot t^{-n} \right) \oplus\CC C$$
and $\CC_c =\CC$ is an  $U(\widehat{\mathfrak{g}}_{(\leqslant 0 )})$-module; the central element $C$ acts on $\CC_c$ as scalar multiplication by $c$  and $\mathfrak{g}\ot t^{-n}$ with $n\leqslant 0$ act trivially. By the results of I. B. Frenkel and Y.-C. Zhu in \cite{FZ} and B.-H. Lian in \cite{Lian}, the space $V_c (\mathfrak{g})$ possesses a vertex algebra structure. Furthermore, any restricted $\widehat{\mathfrak{g}}$-module of level $c$ is naturally a module for the vertex algebra $V_c (\mathfrak{g})$ and vice versa; see \cite[Chapter 6]{LLi} for more details and references. In this paper, we study a certain quantum version of that result for $\mathfrak{g}=\gl_N$.

The notion of {\em quantum vertex operator algebra} was introduced by P. Etingof and D. Kazhdan in \cite{EK}. They constructed examples of quantum vertex operator algebras by quantizing the quasiclassical structure on $V_c(\mathfrak{sl}_N)$ when the classical $r$-matrix on 
$\mathfrak{sl}_N$ is rational, trigonometric or elliptic. The corresponding vertex operator map was defined using quantum current $\Tc$, introduced by N. Yu. Reshetikhin and M. A. Semenov-Tian-Shansky in \cite{RS}, which satisfies a commutation relation of the form\footnote{We explain the precise meaning of \eqref{re_introduction} in Section \ref{sec22}.}
\begin{align}\label{re_introduction}
\Tc_1(u)\R(u-v+hC)^{-1} \Tc_2(v) \R(u-v)=\R(-v+u)^{-1}\Tc_2(v)\R(-v+u-hC)\Tc_1(u);
\end{align}
 cf.   also Ding's realization \cite{D} of the quantum affine algebra in type $A$.
 Later on, the theory of quantum vertex algebras  was further developed and generalized by H.-S. Li; see  \cite{Li_g1,Li05,Li06,Li} and references therein.
Specifically, in  \cite{Li},  an {\em $h$-adic quantum vertex algebra}   was constructed on the level $c$ universal vacuum module  for a certain cover  of the  double Yangian $\wht{\DY(\mathfrak{sl}_2)}$ from \cite{Kho} for any generic $c\in\CC$.  Moreover, it was proved  that any highest weight $\wht{\DY(\mathfrak{sl}_2)}$-module of level $c$ is naturally equipped with a  module structure  for that $h$-adic quantum vertex algebra at the   level $c$.

In this paper, we employ commutation relation \eqref{re_introduction}, where $\R$   denotes the  
normalized Yang $R$-matrix,  $h$ is a formal parameter and $C$  a central element,   to define an associative algebra $\A(\R)$ over the ring $\CC[[h]]$, which we refer to as 
the
{\em quantum current algebra
associated to} $\R$. It is worth noting that the  classical limit of \eqref{re_introduction} coincides with the commutation relation for the affine Lie algebra $\wht{\gl}_N$. 
 We   investigate   properties of the algebra $\A(\R)$ and, in particular, we use the fusion procedure originated in \cite{J} to give explicit  formulae for the elements of its center at the critical level. 

Next,  we  introduce the notion of {\em restricted} $\A(\R)$-module    in parallel with the representation theory of the affine Lie algebras; see, e.g., \cite[Chapter 6]{LLi}. For any complex number $c$ we consider the {\em vacuum module} $\Vc_c(\R)$ of level $c$ for the algebra $\A(\R)$, which presents an example of restricted $\A(\R)$-module.
We show that, as a $\CC[[h]]$-module, $\Vc_c(\R)$   is isomorphic to the $h$-adically completed vacuum module $\Vc_c(\gl_N)$ over the double Yangian for the Lie algebra $\gl_N$.
Hence, due to 
the
aforementioned  Etingof--Kazhdan's construction \cite{EK}, the $\CC[[h]]$-module $\Vc_c(\R)=\Vc_c(\gl_N)$  possesses a  quantum vertex algebra structure. 
This structure was recently studied by N. Jing, A. Molev, F. Yang and the author in \cite{JKMY}, where the center of $\Vc_c(\gl_N)$ was determined, and also in \cite{K}, where a certain connection between quasi $\Vc_c(\gl_N)$-modules and a class of reflection algebras of A. Molev and E. Ragoucy \cite{MR}  was established.

 The main result of this paper, Theorem \ref{main} states that
any module for the quantum vertex algebra $\Vc_c(\R)$ is naturally equipped with a  structure of  restricted $\A(\R)$-module of level $c$ and, conversely, that 
 any    restricted $\A(\R)$-module of level $c$ is naturally equipped with a structure of   module for the quantum vertex algebra $\Vc_c(\R)$. 
Roughly speaking, the proof  of the theorem   relies on the fact that commutation relation \eqref{re_introduction} possesses 
a form similar to
the {\em $\Sc$-locality property}, which is one of the fundamental quantum vertex algebra axioms. 

We should mention that, starting with the work of E. K. Sklyanin \cite{S}, various classes of  {\em reflection algebras}, which are   defined via relations of the form
similar to or same as
\beq\label{reflectionequation}
R(u-v)B_1 (u)R(u+v)B_2 (v)=B_2 (v)R(u+v)B_1 (u)R(u-v),
\eeq
thus resembling commutation relation \eqref{re_introduction},
 were extensively studied. For more details   the reader may consult   \cite{GS,JW,KS,KJC,MR,MRS} and references therein. However, in contrast with 
\eqref{re_introduction}, {\em reflection equation} \eqref{reflectionequation} does not seem to directly give rise to the  $\Sc$-locality property, i.e. to the quantum vertex algebra structure; see \cite{K}.

\section{Preliminaries}\label{sec1}
\numberwithin{equation}{section}

In this section, we
first
 recall some properties of the rational $R$-matrix.  Next, we define    the (completed) double Yangian for the Lie algebra $\gl_N$ and its vacuum module. Finally, we recall the notions of quantum vertex algebra and module for
a quantum vertex algebra, which play a central role in this paper.

\subsection{Rational \texorpdfstring{$R$}{R}-matrix}

\allowdisplaybreaks

Let $N\geqslant 2$ be an integer and let $h$ be a formal parameter.
We follow \cite[Section 2.2]{JKMY} to recall the definition and some basic properties of the rational $R$-matrix over the ring $\CC[[h]]$.
Consider the {\em Yang $R$-matrix} over $\CC[[h]]$, 
\beq\label{yang}
R(u)=1-hPu^{-1}\,\in\, \ndo\CC^N\ot \ndo\CC^N \ts[h,u^{-1}],
\eeq
where $1\colon x\ot y\mapsto x\ot y$ is the identity and $P\colon x\ot y\mapsto y\ot x$ is the permutation operator 
on
 $\CC^N\ot \CC^N $.
There exists a
 unique   series $g(u)$ in $1+u^{-1}\CC[[u^{-1}]]$
such that
\beq\label{newcsym}
g(u+N)=g(u)(1-u^{-2}).
\eeq

The $R$-matrix $\overline{R}(u)=\wvr{R}_{12}(u)=g(u/h)R(u)$ possesses the  {\em unitarity property}
\beq\label{uni}
\wvr{R}_{12}(u)\wvr{R}_{12}(-u)=1.
\eeq
This as well as  Yang $R$-matrix \eqref{yang},   the {\em Yang--Baxter equation}
\beq\label{ybe}
\R_{12}(u)\ts \R_{13}(u+v)\ts \R_{23}(v)
=\R_{23}(v)\ts \R_{13}(u+v)\ts \R_{12}(u).
\eeq
Both sides of \eqref{ybe} are operators on the triple tensor product $(\CC^{N})^{\ot 3}$ and the subscripts indicate the copies of $\CC^N$ on which the $R$-matrices are applied, e.g., 
$\R_{12}(u)=\R(u)\ot 1$.

Due to \eqref{newcsym}, the $R$-matrix  $\R(u)$ possesses the {\em crossing symmetry} properties,
\beq\label{csym}
\left(\wvr{R}_{12}(u)^{-1}\right)^{t_1}\wvr{R}_{12}(u+hN)^{t_1}=1\fand \left(\wvr{R}_{12}(u)^{-1}\right)^{t_2}\wvr{R}_{12}(u+hN)^{t_2}=1,
\eeq
where
$t_i$ denotes the transposition applied on the tensor factor $i=1,2$. As  in \cite[Section 4.2]{JKMY},   we can write \eqref{csym} using the ordered product notation as 
\beq\label{csym_equiv}
\wvr{R}_{12}(u)^{-1}\cdotrl\wvr{R}_{12}(u+hN)=1
\fand
\wvr{R}_{12}(u)^{-1}\cdotlr\wvr{R}_{12}(u+hN)=1,
\eeq
where the subscript RL (LR) in \eqref{csym_equiv} indicates that the first tensor factor of $\wvr{R}_{12}(u)^{-1}$ is applied from the right (left) while the second tensor factor  of $\wvr{R}_{12}(u)^{-1}$ is applied from the left (right).\footnote{Strictly speaking, notation used in \cite{JKMY}  slightly differs. The equalities in \eqref{csym_equiv} are  expressed therein as
$${}^{rl}\wvr{R}_{12}(u)^{-1}\wvr{R}_{12}(u+hN)=1
\fand
{}^{lr}\wvr{R}_{12}(u)^{-1}\wvr{R}_{12}(u+hN)=1.$$ } Indeed, \eqref{csym_equiv} is obtained by applying the transposition $t_1$ on the first and $t_2$ on the second equality in \eqref{csym}.

\subsection{Double Yangian for \texorpdfstring{$\gl_N$}{glN}}\label{sec12}

The {\em double Yangian} $\DY(\gl_N)$ for the Lie algebra $\gl_N$  is  the associative algebra over the ring $\CC[[h]]$ generated by the   central element $C$ and the elements
$t_{ij}^{(\pm r)}$, where $i,j=1,\ldots , N$ and $r=1,2,\ldots$,    subject to the  defining relations  
\begin{align}
R\big(u-v \big)\ts T_1(u)\ts T_2(v)&=T_2(v)\ts T_1(u)\ts
R\big(u-v \big),\label{RTT2}\\
R\big(u-v \big)\ts T^+_1(u)\ts T^+_2(v)&=T^+_2(v)\ts T^+_1(u)\ts
R\big(u-v \big),\label{RTT1}\\
\wvr{R}\big(u-v+hC/2\big)\ts T_1(u)\ts T^+_2(v)&=T^+_2(v)\ts T_1(u)\ts
\wvr{R}\big(u-v-hC/2\big),\label{RTT3}
\end{align}
see \cite{EK4,EK,I,K}.
The elements $T(u),T^+(u)\in \ndo\CC^N \ot \DY(\gl_N)[[u^{\mp 1}]]$ are defined by
\begin{align*}
T(u)=\sum_{i,j=1}^N e_{ij}\ot t_{ij}(u)\Fand T^{+}(u)=\sum_{i,j=1}^N e_{ij}\ot t_{ij}^{+}(u),
\end{align*}
where the $e_{ij}$ denote the matrix units
and the series $t_{ij}(u) $ and $t_{ij}^+ (u) $ are defined by
$$t_{ij}(u)=\delta_{ij}+h\sum_{r=1}^{\infty} t_{ij}^{(r)}u^{-r}\Fand t_{ij}^+ (u)=\delta_{ij}-h\sum_{r=1}^{\infty} t_{ij}^{(-r)}u^{r-1}.$$
We  indicate a copy of the matrix in the tensor product algebra
$ (\ndo\CC^N)^{\ot m}\ot\DY(\gl_N)$ by subscripts,
so that, for example, we have
\beq\label{subscript}
T_k(u)=\sum_{i,j=1}^{N} 1^{\otimes (k-1)} \ot e_{ij} \ot 1^{\ot (m-k)} \ot t_{ij}(u).
\eeq
In particular, we have $m=2$ and $k=1,2$ in defining relations \eqref{RTT2}--\eqref{RTT3}.

The {\em Yangian} $\Y(\gl_N)$ is the subalgebra of  $\DY(\gl_N)$ generated by the elements $t_{ij}^{(r)}$, where $i,j=1,\ldots,N$ and $r=1,2,\ldots$  The {\em dual Yangian } $\Y^+ (\gl_N)$ is the subalgebra of the double Yangian $\DY(\gl_N)$ generated by the elements $t_{ij}^{(-r)}$, where $i,j=1,\ldots,N$ and $r=1,2,\ldots$ 
For any   $c\in\CC$ denote by $\DY_c (\gl_N)$ the {\em double Yangian at the level $c$}, which is defined as the quotient of the algebra $\DY(\gl_N)$ by the ideal generated by the element $C-c$.

For any integer $p\geqslant 1$ let $\I_p(\gl_N)$ be the left ideal in 
$\DY_c(\gl_N)$ generated by  all elements $t_{ij}^{(r)}$, where $i,j=1,\ldots ,N$ and $r\geqslant p$. Introduce the completion of the double Yangian $\DY_c(\gl_N)$ at the level $c$ as the inverse limit
$$\wtld{\DY}_c(\gl_N) \,=\,\lim_{\longleftarrow} \,\DY_c(\gl_N)\,/\, \I_p(\gl_N).$$

\subsection{Vacuum module over the double Yangian}\label{sec13}
Let $V$ be an arbitrary $\CC[[h]]$-module.
 The {\em $h$-adic topology} on  $V$ is the topology generated by the  basis $v+h^n V$, where $v\in V$ and $n\in\mathbb{Z}_{\geqslant 1}$.
Recall that  $V$ is said to be {\em torsion-free} if  $hv\neq 0$ for all nonzero $v\in V$ and that $V$ is said to be {\em separable} if $\cap_{m\geqslant 1} h^m V=0$. A $\CC[[h]]$-module V is said to be {\em topologically free} if it is separable, torsion-free and complete with respect to the $h$-adic topology. For more details on   topologically free $\CC[[h]]$-modules see \cite[Chapter XVI]{Kas}.

We now introduce the vacuum module over the double Yangian as in  \cite[Section 4.2]{JKMY}.
Let $\W_c(\gl_N)$ be  the left ideal in 
$\DY_c(\gl_N)$  generated by the elements $t_{ij}^{(r)}$, where $i,j=1,\ldots ,N$ and $ r=1,2,\ldots$ 
By the Poincar\'e--Birkhoff--Witt theorem for the double Yangian, see \cite[Theorem 2.2]{JKMY}, the quotient
\beq\label{imageofone}
\DY_c(\gl_N) \,/\,\W_c(\gl_N)
\eeq
is isomorphic, as a $\CC[[h]]$-module, to the  dual Yangian $\Y^+(\gl_N)$.
The {\em vacuum module $\Vc_c(\gl_N)$ at the level $c$} over the double Yangian is 
defined as the $h$-adic completion of quotient \eqref{imageofone}.
The  vacuum module $\Vc_c(\gl_N)$ is topologically free   $\wtld{\DY}_c (\gl_N)$-module. We denote by $\vac$ the image of the unit $1\in \DY_c(\gl_N)$ in quotient \eqref{imageofone}.

For positive integers $n$ and $m$ introduce functions depending on the variable $z$ and the
families of variables
$u=(u_1,\dots,u_n)$ and $v=(v_1,\dots,v_m)$ with values in
the space
$(\ndo\mathbb{C}^{N})^{\ot  n} \otimes
(\ndo\mathbb{C}^{N})^{\ot  m}
$
by
\begin{align}
\R_{nm}^{12}(u|v|z)= \prod_{i=1,\dots,n}^{\longrightarrow} 
\prod_{j=n+1,\ldots,n+m}^{\longleftarrow} \R_{ij}(z+u_i -v_{j-n}),\label{rnm12}\\
\cev{\R}_{nm}^{12}(u|v|z)= \prod_{i=1,\dots,n}^{\longleftarrow} 
\prod_{j=n+1,\ldots,n+m}^{\longrightarrow} \R_{ij}(z+u_i -v_{j-n}),\label{rnm123}
\end{align}
where the arrows indicate the order of the factors. For example, we have
$$\R_{22}^{12}(u|v|z)=\R_{14}\R_{13}\R_{24}\R_{23} \fand\cev{\R}_{22}^{12}(u|v|z)=\R_{23}\R_{24}\R_{13}\R_{14},$$
where $\R_{ij}= \R_{ij}(z+u_i -v_{j-n}).$

We adopt the following expansion convention. For any variables $x_1,\ldots , x_k$ expressions of the form $(x_1+\ldots +x_k)^s$ with $s<0$ should be expanded in nonnegative powers of the variables $x_2,\ldots ,x_k$.
In particular, expressions of the form $(z+u_i-v_{j-n})^s$ with $s<0$ should be expanded in  negative powers of the variable $z$, so that \eqref{rnm12} and \eqref{rnm123} contain only nonnegative powers of the variables $u_1,\ldots ,u_n$ and $v_{1},\ldots ,v_m$.
Also, we write
\beq\label{rnm1234}
\R_{nm}^{12}(u|v)=\R_{nm}^{12}(u|v|0)\fand \cev{\R}_{nm}^{12}(u|v)=\cev{\R}_{nm}^{12}(u|v|0),
\eeq
where, due to the aforementioned expansion convention,  expressions of the form $( u_i-v_{j-n})^s$ with $s<0$, which appear  in \eqref{rnm1234}, are    expanded in negative powers of the variable $u_i$, so that they contain only nonnegative powers of  $v_{j-n}$.
The functions $R_{nm}^{12}(u|v|z)$ and $\cev{R}_{nm}^{12}(u|v|z)$
corresponding to Yang $R$-matrix \eqref{yang} can be defined analogously.

Introduce the operators on $(\ndo\CC^N)^{\ot n} \ot \Vc_{c}(\gl_N)$ by
\begin{align*}
T_{[n]}^{+}(u|z)=T_{1}^{+}(z+u_1)\ldots T_{n}^{+}(z+u_n)\fand T_{[n]}(u|z)=T_{1}(z+u_1)\ldots T_{n}(z+u_n).
\end{align*}
Note that, due to our expansion convention, the operator $T_{[n]}(u|z)$ contains only nonnegative powers of the variables $u_1,\ldots ,u_n$.
Also, we  write
\beq\label{rnm1234t}
T_{[n]}^{+}(u )=T_{1}^{+}(u_1)\ldots T_{n}^{+}(u_n)\fand T_{[n]}(u)=T_{1}(u_1)\ldots T_{n}(u_n).
\eeq
Note that both expressions in \eqref{rnm1234t} can be viewed as series with coefficients in the double Yangian.
Using  defining relations  \eqref{RTT2}--\eqref{RTT3}  one can verify the following equations for the operators on 
\beq\label{expl1}
\underbrace{(\ndo\CC^N)^{\ot n}}_{1} \ot \underbrace{(\ndo\CC^N)^{\ot m}}_{2}\ot \underbrace{\Vc_c(\gl_N)}_{3},
\eeq
 which were given in  \cite{EK}:
\begin{align}
&R_{nm}^{12}(u|v|z_1-z_2)T_{[n]}^{+13}(u|z_1)T_{[m]}^{+23}(v|z_2)
=T_{[m]}^{+23}(v|z_2)T_{[n]}^{+13}(u|z_1)R_{nm}^{12}(u|v|z_1-z_2),\label{rtt1}\\ 
&R_{nm}^{12}(u|v|z_1-z_2)T_{[n]}^{13}(u|z_1)T_{[m]}^{23}(v|z_2)
=T_{[m]}^{23}(v|z_2)T_{[n]}^{13}(u|z_1)R_{nm}^{12}(u|v|z_1-z_2),\label{rtt2}\\ 
&\overline{R}_{nm}^{\ts 12}(u|v|z_1-z_2+h c/2)T_{[n]}^{13}(u|z_1)T_{[m]}^{+23}(v|z_2)\nonumber\\
&\qquad\qquad\qquad\qquad\qquad=T_{[m]}^{+23}(v|z_2)T_{[n]}^{13}(u|z_1)\overline{R}_{nm}^{\ts 12}(u|v|z_1-z_2-h c/2).
\label{rtt3}
\end{align}
In relations \eqref{rtt1}--\eqref{rtt3}, we use superscripts to indicate tensor factors  in accordance with
\eqref{expl1}. For example, $T_{[n]}^{+13}(u|z_1)$ is applied on  tensor factors $1,\ldots ,n$ and $n+m+1$, and  $T_{[m]}^{+23}(v|z_2)$ is applied on  tensor factors $n+1,\ldots ,n+m$ and $n+m+1$. We will often use such notation throughout this paper.

\subsection{Quantum vertex algebras}
The definitions in this section are presented in the form which we find to be  suitable for the  setting of this paper. In particular, the next definition of quantum vertex algebra coincides with  \cite[Definition 3.1]{JKMY}. It presents a minor modification of the original definition of quantum vertex operator algebra given by Etingof and Kazhdan in \cite{EK}, as explained in \cite[Remark 3.2 and 3.4]{JKMY}.
For more details on the axiomatics of quantum vertex algebras and related structures the reader may consult \cite{EK,Li}. 
 From now on, the tensor products are understood as $h$-adically completed.
\begin{defn}\label{qvoa}
A {\em quantum vertex algebra} is a quintuple $(V,Y,\vac,D,\Sc)$ which satisfies the following axioms:
\begin{enumerate}[(1)]
\item  $V$ is a topologically free $\mathbb{C}[[h]]$-module.
\item $Y$ is a $\mathbb{C}[[h]]$-module map (the {\em vertex operator map})
\begin{align*}
Y \colon V\ot V&\to V((z))[[h]]\\
u\ot v&\mapsto Y(z)(u\ot v)=Y(u,z)v=\sum_{r\in\mathbb{Z}} u_r v \ts z^{-r-1}
\end{align*}
which satisfies the {\em weak associativity}:
for any $u,v,w\in V$ and $n\in\mathbb{Z}_{\geqslant 0}$
there exists $r\in\mathbb{Z}_{\geqslant 0}$
such that
\begin{equation}\label{associativity}
(z_0 +z_2)^r\ts Y(u,z_0 +z_2)Y(v,z_2)\ts w - (z_0 +z_2)^r\ts Y\big(Y(u,z_0)v,z_2\big)\ts w
\in h^n V[[z_0^{\pm 1},z_2^{\pm 1}]].
\end{equation}
\item $\vac$ is an element of $V$ (the {\em vacuum vector}) which satisfies
\beq\label{v1}
Y(\vac ,z)v=v\quad\text{for all }v\in V,
\eeq
and for any $v\in V$ the series $Y(v,z)\ts\vac$ is a Taylor series in $z$ with
the property
\beq\label{v2}
\lim_{z\to 0} Y(v,z)\ts\vac =v.
\eeq
\item $D$ is a $\mathbb{C}[[h]]$-module map $ V\to V$
 which satisfies
\begin{align}
&D\ts\vac =0\Fand\frac{d}{dz}Y(v,z)=[D,Y(v,z)]\quad\text{for all }v\in V.
\label{d2}
\end{align}
\item $\Sc=\Sc(z)$ is a $\mathbb{C}[[h]]$-module map
$ V\otimes V\to V\otimes V\otimes\mathbb{C}((z))$ which satisfies the {\em shift condition}
\begin{align}
&[D\otimes 1, \mathcal{S}(z)]=-\frac{d}{dz}\mathcal{S}(z),\label{s1}\\
\intertext{the {\em Yang--Baxter equation}}
&\mathcal{S}_{12}(z_1)\ts\mathcal{S}_{13}(z_1+z_2)\ts\mathcal{S}_{23}(z_2)
=\mathcal{S}_{23}(z_2)\ts\mathcal{S}_{13}(z_1+z_2)\ts\mathcal{S}_{12}(z_1),\label{s2}\\
\intertext{the {\em unitarity condition}}
&\mathcal{S}_{21}(z)=\mathcal{S}^{-1}(-z),\label{s3}
\end{align}
 and the $\mathcal{S}$-{\em locality}:
for any $u,v\in V$ and $n\in\mathbb{Z}_{\geqslant 0}$ there exists
$r\in\mathbb{Z}_{\geqslant 0}$ such that  
\begin{align}
&(z_1-z_2)^{r}\ts Y(z_1)\big(1\otimes Y(z_2)\big)\big(\mathcal{S}(z_1 -z_2)(u\otimes v)\otimes w\big)
\nonumber\\
&\quad-(z_1-z_2)^{r}\ts Y(z_2)\big(1\otimes Y(z_1)\big)(v\otimes u\otimes w)
\in h^n V[[z_1^{\pm 1},z_2^{\pm 1}]]\quad\text{for all }w\in V.\label{locality}
\end{align}
\end{enumerate}
\end{defn}

It was proved in \cite{Li} that the {\em $\Sc$-Jacobi identity}
\begin{align}
&z_0^{-1}\delta\left(\frac{z_1 -z_2}{z_0}\right) Y(z_1)(1\ot Y(z_2))(u\ot v\ot w)\non\\
&\qquad-z_0^{-1}\delta\left(\frac{z_2-z_1}{-z_0}\right) Y(z_2)(1\ot Y(z_1)) \left(\Sc(-z_0)(v\ot u)\ot w\right)\non\\
&\qquad\qquad=z_2^{-1}\delta\left(\frac{z_1 -z_0}{z_2}\right)Y(Y(u,z_0)v,z_2 )w\quad\text{for all}\quad u,v,w\in V\label{jacobi}
\end{align}
is equivalent to weak associativity \eqref{associativity} and $\Sc$-locality \eqref{locality}. In what follows, the (appropriately modified)  $\Sc$-Jacobi identity is used to define the notion of  module for a quantum vertex algebra. Originally, modules for $h$-adic nonlocal vertex algebras, which present a generalization of quantum vertex algebras, as well as  modules for   some related structures,  were introduced and studied by Li; see, e.g., \cite{Li_g1,Li05,Li06,Li} and references therein.

\begin{defn}\label{qvoamodule}
 Let $(V,Y,\vac,D,\Sc)$  be a quantum vertex algebra. A {\em   $V$-module} is a pair $(W,Y_W)$, where $W$ is a topologically free $\CC[[h]]$-module
and 
\begin{align*}
Y_W(z)\colon V\ot W&\to W((z))[[h]]\\
v\ot w&\mapsto Y_W(z)(v\ot w)=Y_W(v,z)w=\sum_{r\in\mathbb{Z}} v_r w\ts z^{-r-1}
\end{align*}
is a $\CC[[h]]$-module map
which satisfies the {\em $\Sc$-Jacobi identity}
\begin{align}
&z_0^{-1}\delta\left(\frac{z_1 -z_2}{z_0}\right) Y_W(z_1)(1\ot Y_W(z_2))(u\ot v\ot w)\non\\
&\qquad-z_0^{-1}\delta\left(\frac{z_2-z_1}{-z_0}\right) Y_W(z_2)(1\ot Y_W(z_1)) \left(\Sc(-z_0)(v\ot u)\ot w\right)\non\\
&\qquad\qquad=z_2^{-1}\delta\left(\frac{z_1 -z_0}{z_2}\right)Y_W(Y(u,z_0)v,z_2 )w\quad\text{for all}\quad u,v \in V\text{ and }w\in W,\label{wjacobi}
\end{align}
and  $$Y_W(\vac,z)w=w\quad\text{ for all } w\in W. $$
 Let  $W_1$ be a topologically free $\CC[[h]]$-submodule of $W$. A pair $(W_1,Y_{W_1})$ is said to be a {\em $V$-submodule} of $W$ if $Y_W (v,z)w_1$ belongs to $ W_1$ for all $v\in V$ and $w_1\in W_1$, where $Y_{W_1}$ denotes the restriction and corestriction of $Y_W$,
$$Y_{W_1}(z)=Y_W (z)\big|_{V\ot W_1}^{W_1}   \big. \,\colon\, V\ot W_1\,\to\, W_1 ((z))[[h]].$$
\end{defn}

The next  lemma, which we use in the proof of Theorem \ref{main}, can be viewed as an $h$-adic analogue of \cite[Remark 2.5]{Li06}; see also \cite[Theorem 4.4.5]{LLi}. Roughly speaking, it is a consequence of \cite[Lemma 2.1]{Li_g1} and the fact that for any quantum vertex algebra $V$ and $V$-module $W$  the quotient $V/h^n V$ is a {\em weak quantum vertex algebra} over $\CC$ and $W/h^n W$ is a {\em $V/h^n V$-module} for all $n\geqslant 1$; see \cite{Li} for details.

\begin{lem}\label{usefullemma}
Let $(V,Y,\vac,D,\Sc)$ be a quantum vertex algebra and let  $W$ be a topologically free $\CC[[h]]$-module equipped with $\CC[[h]]$-module map
\begin{align*}
Y_W(z)\colon V\ot W&\to W((z))[[h]]\\
v\ot w&\mapsto Y_W(z)(v\ot w)=Y_W(v,z)w=\sum_{r\in\mathbb{Z}} v_r w z^{-r-1}
\end{align*}
which satisfies $Y_W(\vac,z)w=w$  for all $w\in W$ and the {\em weak associativity}:
for any $u,v\in V$, $w\in W$ and $n\in\mathbb{Z}_{\geqslant 0}$
there exists $r\in\mathbb{Z}_{\geqslant 0}$
such that
\begin{align}
&(z_0 +z_2)^r\ts Y_W (u,z_0 +z_2)Y_W(v,z_2)\ts w\non\\
&\qquad - (z_0 +z_2)^r\ts Y_W\big(Y(u,z_0)v,z_2\big)\ts w
\in h^n W[[z_0^{\pm 1},z_2^{\pm 1}]].\label{associativityw}
\end{align}
Then $(W,Y_W)$ is a $V$-module.
In particular, 
the $\mathcal{S}$-{\em locality} holds, i.e.
for any $u,v\in V$ and $n\in\mathbb{Z}_{\geqslant 0}$ there exists
$r\in\mathbb{Z}_{\geqslant 0}$ such that
\begin{align}
&(z_1-z_2)^{r}\ts Y_W(z_1)\big(1\otimes Y_W(z_2)\big)\big(\mathcal{S}(z_1 -z_2)(u\otimes v)\otimes w\big)
\nonumber\\
&-(z_1-z_2)^{r}\ts Y_W(z_2)\big(1\otimes Y_W(z_1)\big)(v\otimes u\otimes w)
\in h^n W[[z_1^{\pm 1},z_2^{\pm 1}]]\quad\text{for all }w\in W.\label{localityw}
\end{align}
\end{lem}

\begin{prf}
The lemma can be proved by arguing as in \cite[Remark 2.16]{Li} and using  \cite[Remark 2.5]{Li06} and \cite[Proposition 2.24]{Li}.
\end{prf}

\section{Quantum current algebras}\label{sec2}

In this section, we introduce quantum current algebra   and   study its properties.  In particular, we construct an action of the quantum current algebra on the vacuum module over the double Yangian, which is an important ingredient of the proof of our main result, Theorem \ref{main} in Section \ref{sec3}.  In the end, we give explicit formulae for families of
central elements of the quantum current algebra at the critical level.

\subsection{Preliminaries}
Our first goal is to introduce the quantum current algebra $\A(\R)$. 
However, its defining relations, i.e. the coefficients in commutation relation \eqref{re} with respect to the variables $u$ and $v$, are given in terms of certain infinite sums.  In order to handle such expressions and employ commutation relation \eqref{re}, we introduce  the appropriate completion of the corresponding free algebra.

For any integer $N\geqslant 2$ let $\F(N)$ be the associative algebra over the ring $\CC[[h]]$ generated by the elements $1$, $C$ and $\tz_{ij}^{(r)}$, where $i,j=1,\ldots, N$ and $r\in\ZZ$, subject to the following defining relations:
$$C\cdot x=x\cdot C\fand 1\cdot x=x\cdot 1=x\qquad\text{for all }x\in \F(N).$$
Hence $1$ is the unit and $C$ is a central element in  $\F(N)$.
Arrange all $\tz_{i j}^{(r)}$ into  Laurent series
\beq\label{tau}
\tz_{ij}(u)=\delta_{ij}-h\sum_{r\in\ZZ}\tz_{ij}^{(r)}u^{-r-1}\,\in\, \F(N)[[u^{\pm 1}]],\quad\text{where }  i,j=1,\ldots ,N,
\eeq
and introduce the elements $\Tc(u)$ in $\ndo\CC^N\ot\F(N)[[u^{\pm 1}]]$ by
\beq\label{TAU}
\Tc(u)=\sum_{i,j=1}^N e_{ij}\ot \tz_{ij}(u),
\eeq
where the $e_{ij}$ denote the matrix units.

For any integer $n\geqslant 1$ let $\I_n(N)$ be the left ideal in the algebra $\F(N)$ generated by  all elements $\tz_{ij}^{(r)}$, where  $i,j=1,\ldots, N$ and $r\geqslant n-1$. Introduce the completion of the algebra $\F(N)$ as the inverse limit
$$\F' (N) = \lim_{\longleftarrow} \ts\F(N)\ts /\ts \I_n(N).$$
The algebra  $\F' (N)$ is naturally equipped with the $h$-adic topology. 
Let $\Ft(N)$ be the $h$-adic completion of $\F'(N)$, i.e. $\Ft(N)=\F'(N)[[h]]$.
For any integer $n\geqslant 1$ denote by $\I_n^h (N)$ the $h$-adically completed left ideal in $\Ft(N)$  generated by $\I_n(N)$ and $h^n \cdot 1$.
Recall the subscript notation from Section \ref{sec12} and the expansion convention from Section \ref{sec13}. 
Consider the following expressions:
\begin{align*}
&\wvr{\Tc}_{[2]}(u,v)=\Tc_1(u)\R(u-v+hC)^{-1} \Tc_2(v) \R(u-v),\\
&\wndr{\Tc}_{[2]}(u,v)=\R(-v+u)^{-1}\Tc_2(v)\R(-v+u-hC)\Tc_1(u).
\end{align*}

\begin{lem}\label{LR}
The expressions $\wvr{\Tc}_{[2]}(u,v)$ and $\wndr{\Tc}_{[2]}(u,v)$ 
are well-defined elements of
\beq\label{space}
\ndo\CC^N \ot \ndo\CC^N \ot \Ft(N) [[u^{\pm 1},v^{\pm 1}]].
\eeq
Moreover, for any integer $n\geqslant 1$ the elements $\wvr{\Tc}_{[2]}(u,v) $ and $\wndr{\Tc}_{[2]}(v,u)$\footnote{Notice  the    swapped variables in this  term.} modulo $\I_n^h(N) $  belong to
\beq\label{space2}
\ndo\CC^N \ot \ndo\CC^N \ot \F(N) [[u^{\pm 1} ]]((v)).
\eeq
\end{lem}

\begin{prf}
The matrix entries of
$\R(u-v+hC)^{-1}$ and  $\R(u-v)$
belong to $\CC[u^{-1}][[hC,h,v]]$. 
For any integer $n\geqslant 1$  the matrix entries of $\Tc (v)$ modulo $\I_n (N) $ belong to $\F(N)((v))$. Therefore,  the expression $\R(u-v+hC)^{-1} \Tc_2(v) \R(u-v)$ is a well-defined element of \eqref{space} and its matrix entries modulo $\I_n^h(N) $ belong to $ \F(N)[u^{-1}]((v))$.
Hence, by left multiplying this expression by $\Tc_1(u)$, we conclude that   $\wvr{\Tc}_{[2]}(u,v)$ 
is a well-defined element of \eqref{space} and, furthermore, that $\wvr{\Tc}_{[2]}(u,v)$ modulo  $\I_n^h(N) $ belongs to \eqref{space2}.
The corresponding statements for $\wndr{\Tc}_{[2]}(v,u)$   can be verified analogously.
\end{prf}

By Lemma \ref{LR},  there exist elements 
$\wvr{\tz}_{i\ts j\ts k\ts l}^{(r,s)}$ and $\wndr{\tz}_{i\ts j\ts k\ts l}^{(r,s)}$ in  $\Ft(N)$, where $i,j,k,l=1,\ldots ,N$  and  
 $r,s\in \ZZ$, such that
\begin{align*}
&\wvr{\Tc}_{[2]}(u,v)=\sum_{i,j,k,l=1}^N\sum_{r,s\in\mathbb{Z}} e_{ij}\ot e_{kl}\ot \wvr{\tz}_{i\ts j\ts k\ts l}^{(r,s)}\ts u^{-r-1}v^{-s-1},\\
&\wndr{\Tc}_{[2]}(u,v)=\sum_{i,j,k,l=1}^N\sum_{r,s\in\mathbb{Z}} e_{ij}\ot e_{kl}\ot \wndr{\tz}_{i\ts j\ts k\ts l}^{(r,s)}\ts u^{-r-1}v^{-s-1}.
\end{align*}

\subsection{Quantum current algebra}\label{sec22}
We are now prepared to introduce the   quantum current algebra $\A  (\R)$.  
Let $\J(N)$ be the  ideal in the algebra $\Ft(N)$ generated by all elements 
$\wvr{\tz}_{i\ts j\ts k\ts l}^{(r,s)}-\wndr{\tz}_{i\ts j\ts k\ts l}^{(r,s)}$, where $i,j,k,l=1,\ldots ,N$ and $r,s\in\ZZ$. 
Introduce its completion as the inverse limit
$$\J' (N)=\lim_{\longleftarrow} \ts \J(N)\ts /\ts \J(N)\cap \I_n (N).$$
The $h$-adic completion $\wtld{[\J'(N)]}$ of
\beq\label{joten}
[\J'(N)]=\left\{a\in\Ft (N)\,:\, h^n a\in \J'(N)\text{ for some }n\geqslant 0  \right\}
\eeq
  is also an ideal in  $\Ft(N)$.
Define the {\em quantum current algebra}  $\A (\R)$ as the quotient of the algebra $\Ft(N)$ by  the ideal  $\wtld{[\J'(N)]}$,
\beq\label{quotient}
\A (\R)\,=\, \Ft(N)\,/\, \, \wtld{[\J'(N)]}.
\eeq
Clearly, the algebra $\A (\R)$ is equipped with the $h$-adic topology.
\begin{pro}\label{free}
The  algebra $\A (\R)$ is topologically free.
\end{pro}

\begin{prf}
We will prove that the algebra $\A (\R)$ is torsion-free, separated and $h$-adically complete.
Let $a$ be an arbitrary  element in $\A (\R)$. Choose any element $b$ in $\Ft(N)$ such that its image, with respect to the canonical map $\pi\colon\Ft(N)\to \Ft(N)/\wtld{[\J'(N)]}$, equals $a$. 
By \cite[Proposition 3.7]{Li} we have
\beq\label{complid}
\wtld{[\J'(N)]} = \left\{a\in\Ft(N)\,:\, h^n a\in \wtld{[\J'(N)]}\text{ for some }n\geqslant 0  \right\}.
\eeq
Suppose that $ha=0$. Then $hb$ belongs to $\wtld{[\J'(N)]}$, so we conclude by \eqref{complid} that  $b$ belongs to $\wtld{[\J'(N)]}$. This implies $a=\pi(b)=0$, so  the algebra $\A (\R)$ is torsion-free.

Suppose that $a$ belongs to $ h^n \A (\R)$ for all $n\geqslant 1$.  This implies
$$b\,\in\,\wtld{[\J'(N)]} \mod h^n \Ft(N) \qquad\text{for all }n\geqslant 1.$$
Hence we can construct a sequence $(b_n)_n$ in $\wtld{[\J'(N)]} $ such that
$b=\lim_n b_n$, with respect to the $h$-adic topology. Since the ideal  $\wtld{[\J'(N)]} $ is $h$-adically complete, we conclude that $\lim_n b_n=b$ belongs to  $\wtld{[\J'(N)]}  $, which implies $a=\pi(b)=0$. Therefore, $\cap_{n\geqslant 1}h^n \A (\R) =0$, so the algebra $\A (\R)$ is  separated.  

It remains to prove that the algebra $\A (\R)$ is $h$-adically complete. Let $(a_n)_n$ be a Cauchy sequence in $\A (\R)$. There exists an increasing and unbounded sequence of nonnegative integers $(m_n)_n$ such that   $a_{n+1}-a_n$ belongs to $h^{m_n}\A (\R)$ for all $n$.
Let $(b_n)_n$ be any sequence in $\Ft(N)$ such that $\pi(b_n)=a_n$ for all $n$.
There exist elements $c_n\in \Ft(N)$ and $d_n\in \wtld{[\J'(N)]}$ such that $b_{n+1}-b_n = h^{m_n} c_n + d_n$ for all $n$. Consider the sequence $(e_n)_n$ in $\Ft(N)$ defined by 
$e_n = b_n-d_{n-1}-\ldots -d_1$. We have $e_{n+1}-e_{n}= h^{m_n} c_n$ for all $n$, so   $(e_n)_n$ is a Cauchy sequence in $\Ft(N)$. Since $\Ft(N)$ is $h$-adically complete, there exists $e\in \Ft(N)$ such that $\lim_n e_n =e$. Hence, for every $k\geqslant 0$ there exists $l\geqslant 1$ such that  $e-e_n\in h^k \Ft(N)$ for all $n\geqslant l$. 
Since $\pi(e_n)=\pi(b_n)=a_n$, this implies $\pi(e)-a_n\in h^k \A (\R)$
for all $n\geqslant l$. Therefore, the limit of the Cauchy sequence $(a_n)_n$ does exist in $\A (\R)$ (and  is  equal to $\pi(e)$), so we conclude that the algebra $\A (\R)$ is $h$-adically complete.
\end{prf}

From now on, we denote the images of the elements $1$, $C$ and $\tz_{ij}^{(r)}$ 
in quotient \eqref{quotient} by $1$, $C$ and $\tz_{ij}^{(r)}$  respectively. 
Also, we  denote by $\tz_{ij}(u)$ and $\Tc(u)$  the corresponding series 
in $\A(\R)[[u^{\pm 1}]]$ and   $\ndo\CC^N \ot \A(\R) [[u^{\pm 1}]]$. 
Defining relations for the   quantum current algebra,
$$\wvr{\tz}_{i\ts j\ts k\ts l}^{(r,s)}-\wndr{\tz}_{i\ts j\ts k\ts l}^{(r,s)}=0,\qquad \text{where}\qquad i,j,k,l=1,\ldots ,N\fand r,s\in\ZZ,$$ 
 can be expressed  as a
commutation relation in $\ndo\CC^N \ot \ndo\CC^N \ot \A (\R) [[u^{\pm 1},v^{\pm 1}]]$,
\begin{align}\label{re}
\Tc_1(u)\R(u-v+hC)^{-1} \Tc_2(v) \R(u-v)=\R(-v+u)^{-1}\Tc_2(v)\R(-v+u-hC)\Tc_1(u).
\end{align}
Since the images of the elements $\wvr{\tz}_{i\ts j\ts k\ts l}^{(r,s)}$ and $\wndr{\tz}_{i\ts j\ts k\ts l}^{(r,s)}$ in quotient \eqref{quotient} coincide,   we   denote them by $\tz_{i\ts j\ts k\ts l}^{(r,s)}$. We also  write $\Tc_{[1]}(u)=\Tc (u)$ and
\beq\label{tedva}
\Tc_{[2]}(u,v)=\sum_{i,j,k,l=1}^N\sum_{r,s\in\mathbb{Z}} e_{ij}\ot e_{kl}\ot  \tz_{i\ts j\ts k\ts l}^{(r,s)}\ts u^{-r-1}v^{-s-1}\,\in\,  (\ndo\CC^N)^{\ot 2} \ot \A(\R) [[u^{\pm 1},v^{\pm 1}]],
\eeq
so that, in particular,  the both sides of commutation relation \eqref{re} coincide with  $\Tc_{[2]}(u,v)$.

Denote by $\I_n^h (\R)$ and $\I_n (\R)$ the images of   the left ideals $\I_n^h (N)  $ and $\I_n (N)  $ in the algebra $\A(\R)$ with respect to the canonical map $\Ft(N) \to \A(\R)$.
\begin{kor}\label{LRkor}
For any $n\geqslant 1$ the element $\Tc_{[2]}(u,v)$  modulo $\I_n^h (\R)$
 belongs to
$$\ndo\CC^N \ot \ndo\CC^N \ot\A(\R)((u,v)).$$
\end{kor}

\begin{prf}
By Lemma \ref{LR}, the  left hand side of  commutation relation \eqref{re} modulo $\I_n^h (\R)$ belongs to
$(\ndo\CC^N)^{\ot 2}\ot\A(\R)[[u^{\pm 1}]]((v))$ while the  right hand side of \eqref{re} modulo $\I_n^h (\R)$ belongs to
$(\ndo\CC^N)^{\ot 2}\ot\A(\R)[[v^{\pm 1}]]((u))$.  Hence, both sides of \eqref{re} possess only finitely many negative powers of the variables $u$ and $v$ modulo  $\I_n^h (\R)$, so the corollary  follows.
\end{prf}

For any integer $n\geqslant 2$ and complex number $a $   introduce the functions with values in $(\ndo\CC^N)^{\ot n}$ in the variables $u=(u_1,\ldots,u_n)$ by
\begin{align}
&\R_{[n,a]}(u)=\prod_{i=1,\ldots ,n-1}^{\longrightarrow}\prod_{j=i+1,\ldots ,n}^{\longrightarrow}\R_{ij}(u_i-u_j+ah),\label{rovi1}\\
&\cev{\R}_{[n,a]}(u)=\prod_{i=1,\ldots ,n-1}^{\longleftarrow}\prod_{j=i+1,\ldots ,n}^{\longleftarrow}\R_{ij}(u_i-u_j+ah),\label{rovi2}
\end{align}
where the arrows indicate the order of the factors.
If $a=0$, we omit the second subscript and write $\R_{[n]}(u)=\R_{[n,0]}(u)$ and $\cev{\R}_{[n]}(u)=\cev{\R}_{[n,0]}(u)$. The functions 
$R_{[n,a]}(u)$ and $\cev{R}_{[n,a]}(u)$ associated with Yang $R$-matrix \eqref{yang} can be defined analogously.
Introduce the expression
\beq\label{tntn67}\Tc_{[n]}(u)=\hspace{-4pt}\prod_{i=1,\ldots ,n}^{\longrightarrow}\hspace{-4pt} \left(\Tc_{i}(u_i)\R_{i\ts i+1}(u_i -u_{i+1}+hC)^{-1}\ldots \R_{i\ts n}(u_i -u_{n}+hC)^{-1}  \right)\,\cdot\, \cev{\R}_{[n]}(u).
\eeq
For example, by setting $\wht{R}_{ij}=\R_{i j}(u_i -u_{j}+hC)^{-1}$, $\R_{i j}=\R_{ij}(u_i-u_j)$ and $n=4$,  we get
$$\Tc_{[4]}(u)=\Tc_{1}(u_1)\wht{R}_{12}\wht{R}_{13}\wht{R}_{14}\Tc_{2}(u_2)\wht{R}_{23}\wht{R}_{24}\Tc_{3}(u_3)\wht{R}_{34}\Tc_{4}(u_4)\R_{34}\R_{24}\R_{23}\R_{14}\R_{13}\R_{12}.$$

\begin{pro}\label{teemlema}
For any  $m\geqslant 1$ the expression $\Tc_{[m]}(u)$ is a well-defined element of 
$$(\ndo\CC^N)^{\ot m} \ot \A(\R)[[u_{1}^{\pm 1},\ldots ,u_{m}^{\pm 1}]].$$ 
Moreover, for any integer $n\geqslant 1$  the element $\Tc_{[m]}(u)$ modulo $\I_n^h (\R)$
 belongs to
$$(\ndo\CC^N)^{\ot m} \ot \A(\R)((u_1,\ldots,u_m)).$$
\end{pro}

\begin{prf}
Fix an integer $n\geqslant 1$.
The proposition follows by induction on $m$. The case $m=1$ is clear (and the case $m=2$ is already given by   Corollary \ref{LRkor}). 
Assume that the proposition holds for some integer $m\geqslant 1$. 
Definition in \eqref{tntn67} implies
\beq\label{tntn68}
\Tc_{[m+1]}(u,v)=\Tc^{13}_{[1]}(u) \,\R_{1m}^{12}(u+hC|v)^{-1} \,\Tc_{[m]}^{23}(v)\,\R_{1m}^{12}(u|v),
\eeq
where $v=(v_1,\ldots, v_m)$, $u$ is a single variable and, as explained in Section \ref{sec13}, the superscripts indicate tensor copies in
$$\underbrace{\ndo\CC^N}_{1} \ot \underbrace{(\ndo\CC^N)^{\ot m}}_{2} \ot \underbrace{\A(\R)}_{3}.$$
By arguing as in the proof of Lemma \ref{LR}, one can prove that
\eqref{tntn68}
is a well-defined element of 
$$\ndo\CC^N \ot (\ndo\CC^N)^{\ot m} \ot \A(\R)[[u^{\pm 1},v_{1}^{\pm 1},\ldots ,v_{m}^{\pm 1}]].$$ 

As for the second assertion, by using    Yang--Baxter equation \eqref{ybe} and commutation relation \eqref{re} one can prove the equality
\begin{align}
&\Tc^{13}_{[1]}(u) \,\R_{1m}^{12}(u+hC|v)^{-1} \,\Tc_{[m]}^{23}(v)\,\R_{1m}^{12}(u|v)\non\\
&\qquad=\R_{1m}^{12\ts r}(u|v)^{-1}\,\Tc_{[m]}^{23}(v) \,\R_{1m}^{12\ts r}(u-hC|v)\,\Tc^{13}_{[1]}(u).\label{sunday}
\end{align}
The superscript $r$ in \eqref{sunday} indicates that the rational functions $\R_{1m}^{12\ts r}(u|v)^{-1}=\R_{1m}^{12}(u|v)^{-1}$
and   $\R_{1m}^{12\ts r}(u-hC|v)=\R_{1m}^{12}(u-hC|v)$  should be expanded in nonnegative powers of the variable $u$, thus violating the expansion convention introduced in Section \ref{sec13}. More precisely, in terms of the aforementioned expansion convention, we have
\begin{align*}
&\R_{1m}^{12\ts r}(u|v)^{-1}=\R_{12}(-v_1+u)^{-1} \ldots \R_{1\ts m+1}(-v_m+u)^{-1}\fand\\
&\R_{1m}^{12\ts r}(u-hC|v)=\R_{1\ts m+1}(-v_m+u-hC)\ldots \R_{12}(-v_1+u-hC),
\end{align*}
while, by  definition in  \eqref{rnm12}, we have
\begin{align*}
&\R_{1m}^{12}(u|v)^{-1}=\R_{12}(u-v_1)^{-1} \ldots \R_{1\ts m+1}(u-v_m)^{-1}\fand\\
&\R_{1m}^{12}(u-hC|v)=\R_{1\ts m+1}(u-v_m-hC)\ldots \R_{12}(u-v_1-hC).
\end{align*}
Since the left hand side of \eqref{sunday} coincides with
$\Tc_{[m+1]}(u,v)$, we may now proceed as in the proof of Corollary \ref{LRkor}. More precisely, the right hand side of \eqref{sunday} modulo  $\I_n^h (\R)$ possesses finitely many negative powers of the variable $u$. By the induction hypothesis, the left hand side of \eqref{sunday} modulo  $\I_n^h (\R)$ possesses finitely many negative powers of the variables $v_1,\ldots ,v_m$. Hence, $\Tc_{[m+1]}(u,v)$ modulo  $\I_n^h (\R)$ belongs to
$$\ndo\CC^N \ot (\ndo\CC^N)^{\ot m} \ot \A(\R)((u,v_1,\ldots,v_m)),$$
as required.
\end{prf}

By Proposition \ref{teemlema}, for any integer $n\geqslant 1$ there exist elements 
$\tz_{i_1\ts j_1\ldots i_n\ts j_n}^{(r_1,...,r_n)}$  in  $\A(\R)$, where $
i_1,j_1,\ldots ,i_n,j_n=1,\ldots ,N$ and  $r_1,\ldots , r_n\in\ZZ$, 
such that
\beq\label{teen}
\Tc_{[n]}(u)=\sum_{i_1,j_1,\ldots,i_n,j_n=1}^N\sum_{r_1,\ldots ,r_n\in\mathbb{Z}} e_{i_1 j_1}\ot\ldots \ot e_{i_n j_n}\ot  \tz_{i_1\ts j_1\ldots i_n\ts j_n}^{(r_1,...,r_n)}\ts u^{-r_1-1}_1\ldots u_n^{-r_n-1}.
\eeq
For any $a\in h\CC[h,C]$ and the variables $u=(u_1,\ldots ,u_n)$ denote by $u+a$ the shifted variables $(u_1 +a,\ldots, u_n +a)$.
We now employ elements \eqref{teen} to write the more general form of commutation relation \eqref{re} in
$$(\ndo\CC^N)^{\ot n} \ot (\ndo\CC^N)^{\ot m} \ot \A(\R)[[u_{1}^{\pm 1},\ldots ,u_{n}^{\pm 1},v_{1}^{\pm 1},\ldots ,v_{m}^{\pm 1}]].$$
\begin{pro}\label{reflectionprop}
For any integers $n,m\geqslant 1$ and the variables $u=(u_1,\ldots ,u_n)$ and $v=(v_1,\ldots ,v_m)$  we have
\begin{align}
\Tc_{[n+m]}(u,v)&=\Tc_{[n]}^{13}(u) \,\R_{nm}^{12}(u+hC|v)^{-1} \,\Tc_{[m]}^{23}(v)\,\R_{nm}^{12}(u|v),\label{re1}
\\
\Tc_{[n+m]}(u,v)&=\R_{nm}^{12\ts r}(u|v)^{-1}\,\Tc_{[m]}^{23}(v) \,\R_{nm}^{12\ts r}(u-hC|v)\,\Tc_{[n]}^{13}(u).
\label{re2}
\end{align}
The superscript $r$ in \eqref{re2} indicates that the rational functions $\R_{nm}^{12\ts r}(u|v)^{-1}=\R_{nm}^{12}(u|v)^{-1}$ and  $\R_{nm}^{12\ts r}(u-hC|v)=\R_{nm}^{12}(u-hC|v)$  are expanded in nonnegative powers of the variables $u=(u_1,\ldots ,u_n)$.
\end{pro}

\begin{prf}
It is clear from definition in \eqref{tntn67} that \eqref{re1}  holds for any integers $n,m\geqslant 1$. Moreover, \eqref{re2} holds for $n=1$; recall  \eqref{sunday}. In general, for  a fixed integer $m\geqslant1$,    one can  prove by induction on $n$,  which relies on   \eqref{sunday} and   Yang--Baxter equation \eqref{ybe}, that the right hand sides of \eqref{re1} and \eqref{re2} coincide.
\end{prf}

\subsection{Quantum current algebra at the level \texorpdfstring{$c$}{c} }\label{sec2222}

For any $c\in \CC$ define the {\em  quantum current algebra  at the level $c$} as the quotient $\A_c(\R)$  of the algebra $\A(\R)$   by the  ideal generated by $C-c$. We now use the vacuum module over the  double Yangian at the level $c$, as defined in Section \ref{sec13}, to obtain an example of an $\A_c(\R)$-module.

\begin{pro}
For any $c\in\CC$  the assignments
\beq\label{homohomo}
\Tc_{[n]}(u)\, \mapsto\, T_{[n]}^+ (u)T_{[n]} (u+hc/2)^{-1}
\eeq
with $n\geqslant 1$ and the variables $u=(u_1,\ldots ,u_n)$
define  a structure of $\A_c(\R)$-module on the vacuum module $\Vc_c(\gl_N)$. In particular, the assignments
\beq\label{chhomoepi}
\Tc_{[n]}(u)\, \mapsto\, T_{[n]}^+ (u)\vac
\eeq
with $n\geqslant 1$  
define a $\CC[[h]]$-module epimorphism 
\begin{align}
\A_c(\R)\,\to\, \Vc_{c}(\gl_N)
.\label{chhomo}
\end{align}
\end{pro}

\begin{prf}
Let us prove that the operator  $T^+(u) T(u+hc/2)^{-1}$ on  
$ \Vc_c(\gl_N)$
 satisfies commutation equation \eqref{re} at the level $c$. Consider the expression
\beq\label{eqq1}
T^+_1(u) T_1(u+hc/2)^{-1}\R(u-v+hc)^{-1} T_2^+(v) T_2(v+hc/2)^{-1}\R(u-v),
\eeq
which corresponds to the left hand side of \eqref{re}. Due to \eqref{RTT3}, this equals to  
$$T^+_1(u)  T_2^+(v)\R(u-v)^{-1} T_1(u+hc/2)^{-1} T_2(v+hc/2)^{-1}\R(u-v).$$
Next, by using \eqref{RTT2} and then canceling the $R$-matrices $\R(u-v)^{\pm 1}$ we get
\beq\label{eqq2}
T^+_1(u)  T_2^+(v)  T_2(v+hc/2)^{-1}T_1(u+hc/2)^{-1}.
\eeq

Consider the expression
\beq\label{eqq3}
\R(-v+u)^{-1}T_2^+(v) T_2(v+hc/2)^{-1}\R(-v+u-hc)T^+_1(u) T_1(u+hc/2)^{-1},
\eeq
which corresponds to the right hand side of \eqref{re}. Due to unitarity property \eqref{uni}, the $R$-matrix $\R(-v+u-hc)$ is equal to $\R(v-u+hc)^{-1} $, so we can employ \eqref{RTT3} and write the given expression as
$$\R(-v+u)^{-1}T_2^+(v)T^+_1(u) \R(v-u)^{-1} T_2(v+hc/2)^{-1} T_1(u+hc/2)^{-1}.$$
Finally, by using \eqref{RTT1} and then canceling the $R$-matrices  $\R(v-u)^{-1}=\R(-v+u)$ and $\R(-v+u)^{-1}$ we obtain \eqref{eqq2}.

Since \eqref{eqq1} and \eqref{eqq3} coincide,  the assignment  $\Tc(u)\mapsto T^+(u) T(u+hc/2)^{-1}$ defines  an  $\A_c(\R)$-module structure  on  $\Vc_c(\gl_N)$. Our previous calculation  also shows that the action of $\Tc_{[2]}(u)$
on $\Vc_c(\gl_N)$ is given by $T_{[2]}^+ (u)T_{[2]} (u+hc/2)^{-1}$, i.e. that
\eqref{homohomo} holds for $n= 2$ as well. In general,  \eqref{homohomo} with $n>2$ can be verified by induction on $n$, which relies on defining relations \eqref{RTT2}--\eqref{RTT3} for the double Yangian. 

As for the second statement, observe that  $T(u)\vac=\vac$, so the composition of 
 map \eqref{homohomo} 
and
the  evaluation map  
$\ndo \Vc_{c}(\gl_N) \ni a \mapsto  a\cdot \vac\in \Vc_{c}(\gl_N)$
coincides with $\CC[[h]]$-module map \eqref{chhomo}.
Finally, surjectivity follows from
the Poincar\'e--Birkhoff--Witt theorem for the double Yangian; see  \cite[Theorem 2.2]{JKMY}.
\end{prf}

\begin{rem}
Action \eqref{homohomo} of the algebra $\A_c(\R)$   possesses the  form of the {\em quantum currents} from \cite{RS}; see also \cite{I, EK}. It suggests, together with some other results of the paper,  a possibility of establishing a closer connection between the quantum current algebra and the (completed) double Yangian for the Lie algebra $\gl_N$.
\end{rem}

Let $A_1^c\subset \A_c(\R)$ be the set which contains the unit $1$ and all monomials 
$$\tz_{i_1\ts j_1}^{(r_1)}  \ldots \tz_{i_n\ts j_n}^{(r_n)},\quad\text{where}\quad 
 r_1,\ldots , r_n\in\ZZ,\, i_1,\ldots ,i_n,j_1,\ldots ,j_n=1,\ldots ,N,\, n=1,2,\ldots $$
Next, let $A_2^c\subset \A_c(\R)$ be the set which contains the unit $1$ and all elements 
$$\tz_{i_1\ts j_1\ts\ldots\ts  i_n\ts j_n }^{(r_1,\ldots ,r_n)},\quad\text{where}\quad 
 r_1,\ldots , r_n\in\ZZ,\, i_1,\ldots ,i_n,j_1,\ldots ,j_n=1,\ldots ,N,\, n=1,2,\ldots $$
 Throughout this paper, $\xpan A$  will  always denote the linear span of a subset $A$ of some  $\CC[[h]]$-module, with respect to the ring $\CC[[h]]$. For any subset $A$ of the algebra $\A_c(\R)$  set
$$[A]=\left\{a\in\A_c(\R)\,:\,h^n a\in A\text{ for some }n\geqslant 0\right\}.$$
We denote the image of the left ideal $\I_p^h (\R)\subset\A(\R)$ with respect to the canonical map   $\A(\R) \to \A_c(\R)$ by $\I_p^h (\R)$

\begin{pro}\label{pro_1}
Let $i=1,2$. For any integer $p\geqslant 1$ and element  $a\in\A_c(\R)$  there exist   $a_i\in [\xpan A_i^c]$ such that $a-a_i$ belongs to $\I_p^h(\R) $.
\end{pro}

\begin{prf}
Fix an integer $p\geqslant 1$.
Let $a$ be an arbitrary element of $ \A_c(\R)$.
Denote by $B_1$  the subset of $\F(N)$   which contains the unit $1\in \F(N)$ and all monomials 
$$\tz_{i_1\ts j_1}^{(r_1)}\ldots \tz_{i_n\ts j_n}^{(r_n)}\in\F(N),\text{ where } 
 r_1,\ldots , r_n\in\ZZ,\, i_1,\ldots ,i_n,j_1,\ldots ,j_n=1,\ldots ,N,\, n=1,2,\ldots$$
Choose any element $b$ in $\Ft (N)$ such that its image in the algebra $\A_c(\R)$, with respect to the canonical map $\pi_c\colon\Ft (N)\to \A_c(\R)$, is equal to $a$. 
There exist an element $b_1$  in 
$$[\xpan B_1]=\left\{x\in\Ft (N)\,:\,h^n x\in \textstyle\xpan B_1\text{ for some }n\geqslant 0\right\}$$
such that
$b-b_1$ belongs to $\I_p^h (N)$. Clearly, its image $a_1=\pi_c(b_1 )$  belongs to $[\xpan A_1^c]$ and satisfies $a-a_1 \in \I_p^h (\R)$.

For any $n\geqslant 1$ and the variables $u=(u_1,\ldots ,u_n)$ the element $\Tc_{[n]}(u)$ can be written as
\begin{align}\label{teem2}
\Tc_{[n]}(u)=\left(\left(\cev{\R}_{[n,C]}(u)\right)^{-1}\cdotrl
 \left(\Tc_{1}(u_1)\ldots \Tc_{n}(u_n)  \right)\right)\cdot \cev{\R}_{[n]}(u).
\end{align}
By using unitarity property \eqref{uni} and crossing symmetry property \eqref{csym_equiv} we can  move all $R$-matrices in \eqref{teem2} to the left hand side, thus getting
\begin{align}\label{teem3}
\cev{\R}_{[n,C+N]}(u) \cdotrl\left(\Tc_{[n]}(u)\cdot \left(\cev{\R}_{[n]}(u)\right)^{-1}\right)=\Tc_{1}(u_1)\ldots \Tc_{n}(u_n)  .
\end{align}
The elements of the set $A_1^c$ are exactly the coefficients, with respect to the variables $u_1,\ldots ,u_n$ (and parameter $h$),  of the matrix entries of $\Tc_{1}(u_1)\ldots \Tc_{n}(u_n)$ while the elements of the set $A_2^c$ are exactly the coefficients, with respect to the variables $u_1,\ldots ,u_n$,  of the matrix entries of $\Tc_{[n]}(u)$. Therefore,  \eqref{teem3} implies that there exist $a_2\in [\xpan A_2^c]$ such that $a_1-a_2$ belongs to $\I_p^h (\R)$. Finally, we obtain
$$a-a_2=(a-a_1)+(a_1 -a_2) \in \I_p^h (\R), $$
as required.
\end{prf}

For the variables  $u=(u_1,\ldots,u_n)$ and an integer $i=1,\ldots,n-1$ set
$$u_{i\leftrightarrow i+1}=(u_1,\ldots,u_{i-1},u_{i+1},u_i,u_{i+2},\ldots,u_n)\fand\cev{u}=(u_n,\ldots ,u_1).$$
We use the following lemma  in the proof of Theorem \ref{izomorfizam}.
\begin{lem}\label{newl}
For any integer $n\geqslant 2$ and the variables  $u=(u_1,\ldots,u_n)$ we have
\beq\label{new}
\R_{i\ts i+1}(u_i -u_{i+1}) \Tc_{[n]}(u) \R_{i\ts i+1}(u_i -u_{i+1})^{-1}
=P_{i\ts i+1} \Tc_{[n]}(u_{i\leftrightarrow i+1})    P_{i\ts i+1}.
\eeq
In particular, for $n=2$  we have $u=(u_1,u_2)$  and
\beq\label{new2}
\R (u_1 -u_{2}) \Tc_{[2]}(u_1,u_2) \R (u_1 -u_{2})^{-1}
=P \ts \Tc_{[2]}(u_{2},u_1)\ts    P .
\eeq
\end{lem}

\begin{prf}
Equality \eqref{new2}  follows    from unitarity property \eqref{uni} and   relation \eqref{re}. Finally, by using \eqref{new2} and Yang--Baxter equation \eqref{ybe}  one can easily  verify \eqref{new}.
\end{prf}

\begin{rem}
One can show by a  short calculation that commutation relation \eqref{re} can be written as $h^2  X=0$, where $X$ belongs to $(\ndo \CC^N)^{\ot 2} \ot \A_c(\R)[[u^{\pm 1},v^{\pm 1}]]$. Due to Proposition \ref{free},  the algebra $\A_c(\R)$ is torsion-free, so we have $X=0$. By considering the coefficients of the matrix entries $e_{ij}\ot e_{kl}$ with respect to the variables $u^{-r-1}v^{-s-1}$ in the equality $X\left|_{h=0}\right.=0$, where $X\left|_{h=0}\right.$ denotes the evaluation of $X$ at $h=0$,  we find
$$[\tz_{ij}^{(r)},\tz_{kl}^{(s)}]=\delta_{jk}\tz_{il}^{(r+s)}-\delta_{il}\tz_{kj}^{(r+s)} +\delta_{r+s\ts 0}\ts c\ts r\left(  \delta_{il}\delta_{jk}-\frac{\delta_{ij}\delta_{kl}}{N}\right),$$
i.e.
the commutation relation for the affine Lie algebra $\wht{\gl}_N =\gl_N\ot\CC[t,t^{-1}]\oplus \CC C$ at the level $C=c$.
\end{rem}

\subsection{Central elements of the quantum current algebra at the critical level}
This section presents a   digression from the main topic  as its results are not needed in the rest of the paper.
Consider the following  permutation operator on  $( \CC^N)^{\ot n}$,  
$$P_{[n]} \colon x_1\ot\ldots \ot x_n\mapsto x_n\ot\ldots \ot x_1.$$
Denote the element $P_{[n]}\ts \Tc_{[n]}(\cev{u})\ts P_{[n]}$ more briefly by  $\cev{\Tc}_{[n]}(u)$, so that, in particular,   \eqref{new2} can be written in the form which resembles $RTT$ relations, 
\beq\label{new22}
\R (u_1 -u_{2}) \Tc_{[2]}(u)  
=\cev{\Tc}_{[2]}(u) \R (u_1 -u_{2}).
\eeq
Recall notation \eqref{rovi1}--\eqref{rovi2}.
 By combining Lemma \ref{newl} and Yang--Baxter equation \eqref{ybe}  one can generalize \eqref{new22} as follows.
\begin{lem}\label{fusionkor}
For any integer $n\geqslant 2$ and the variables  $u=(u_1,\ldots,u_n)$ we have 
\beq\label{fusionrtz}
R_{[n]}(u)\ts \Tc_{[n]}(u)\ts  \ts=\ts\cev{\Tc}_{[n]}(u) R_{[n]}(u).
\eeq
\end{lem}

We now recall a special version of  the fusion procedure for  Yang $R$-matrix \eqref{yang} originated in \cite{J}; see also \cite[Section 6.4]{M} for more details and references.
The symmetric group $\Sym_n$ acts on the space $(\CC^N)^{\ot n}$ by permuting the tensor factors.
Denote by $A^{(n)}$  the action  of the anti-symmetrizer  
$$
a^{(n)}=\frac{1}{n!}\sum_{s\in\Sym_n} \sgn \ts s\cdot s\,\in\,\CC[\Sym_n]
$$
on   $(\CC^N)^{\ot n}$. Due to \cite{J}, the consecutive evaluations $u_1 =0, u_2=-h, u_3=-2h, \ldots , u_n =-(n-1)h$ 
of the variables $u=(u_1,\ldots ,u_n)$ in $R_{[n]}(u)$
 are well-defined and 
   we have
\beq\label{fusion}
R_{[n]}(u)\big|_{u_1=0}
\big|_{u_2=-h}\ldots \big|_{u_n=-(n-1)h}=n! \ts A^{(n)}.
\eeq

For the variable $u$   set
$u_{[n]} =(u,u-h,\ldots ,u -(n-1) h).$
 The next two lemmas are used in the proof of Theorem \ref{centralthm}.
\begin{lem}\label{superrr}
The following equalities hold on $\ndo\CC^N \ot (\ndo\CC^N)^{\ot n}$:
\begin{align} 
A^{(n)}\ts \wvr{R}_{1n}^{12}(u_0|u_{[n]})&= \cev{\wvr{R}}_{1n}^{12}(u_0|u_{[n]})\ts A^{(n)},
\label{xxx1}\\
A^{(n)} \ts \wvr{R}_{1n}^{12\ts r}(u_0|u_{[n]})&= \cev{\wvr{R}}_{1n}^{12\ts r}(u_0|u_{[n]})\ts A^{(n)},
\label{xxx2}
\end{align}
where 
the tensor factors of $\ndo\CC^N \ot (\ndo\CC^N)^{\ot n}$ are labeled by
$0,1,\ldots ,n$ and 
  $A^{(n)}$  is applied on the tensor factors $1,\ldots,n$, i.e. $A^{(n)}$ denotes the operator $1\ot A^{(n)}$ on $\ndo\CC^N \ot (\ndo\CC^N)^{\ot n}$.
	The superscript $r$ in  \eqref{xxx2} indicates the rational functions $\wvr{R}_{1n}^{12\ts r}(u_0|u_{[n]})=\wvr{R}_{1n}^{12\ts }(u_0|u_{[n]})$
and   $ \cev{\wvr{R}}_{1n}^{12\ts r}(u_0|u_{[n]})= \cev{\wvr{R}}_{1n}^{12}(u_0|u_{[n]})$  are expanded in negative powers of the variable  $u$.
\end{lem}

\begin{prf}
The lemma follows directly from fusion procedure \eqref{fusion}. More details on its proof can be found in, e.g., \cite[Lemma 3.1]{K}.
\end{prf}

\begin{lem}
The following equality holds on $  (\ndo\CC^N)^{\ot n}\ot \A_{-N}(\R)$:
\beq\label{xxx3}
A^{(n)} \ts \Tc_{[n]} (u_{[n]})\ts=\ts   \cev{\Tc}_{[n]} (u_{[n]})\ts A^{(n)}.
\eeq
\end{lem}

\begin{prf}
The   lemma easily follows by applying the consecutive evaluations $u_1 =0, u_2 =-h, \ldots , u_n =-(n-1)h$ on \eqref{fusionrtz} and  employing fusion procedure \eqref{fusion}.
\end{prf}

We now consider the quantum current algebra   at the critical level $c=-N$.
For each $n=1,\ldots , N$ introduce the series
\beq\label{center0}
\Tf_{[n]} (u)\,=\,\tr_{1,\ldots ,n} \ts A^{(n)}\ts \Tc_{[n]}(u_{[n]})\,\in\,\A_{-N}(\R)[[u^{\pm 1}]],
\eeq
where the trace is taken over all $n$   copies of $\ndo \CC^N$. 
The proof of the next theorem is   similar  to the proof of \cite[Theorem 3.2]{K}. It relies on certain techniques whose $RTT$ counterparts are  well-known; see, e.g., proof of \cite[Theorem 3.2]{FJMR} or \cite[Theorem 4.4]{JKMY}, cf. also \cite{O,T}.

\begin{thm}\label{centralthm}
All coefficients of $\Tf_{[n]} (u)$ belong to the center of the algebra $\A_{-N}(\R)$.
\end{thm}

\begin{prf}
We will prove the equality 
\beq\label{center2}
\Tc (u_0)\ts\Tf_{[n]} (u)\ts=\ts\Tf_{[n]} (u)\ts\Tc (u_0) 
\eeq
in $\ndo\CC^N \ot \A_{-N}(\R)$,
which  implies the statement of the theorem. By applying
$\Tc (u_0)$ on \eqref{center0} we get
\beq\label{xx1}
\tr_{1,\ldots ,n} \ts A^{(n)}\ts \Tc_0 (u_0) \Tc_{[n]}(u_{[n]}),
\eeq
where the expression under the trace belongs to $\ndo\CC^N \ot (\ndo\CC^N)^{\ot n} \ot \A_{-N}(\R)$.
The copies of $\ndo\CC^N$ in \eqref{xx1} are labeled by $0,1,\ldots ,n$. The series $\Tc (u_0)$ is applied on the tensor factor $0$ while  $A^{(n)}$ and $  \Tc_{[n]}(u_{[n]})$ are applied on the tensor factors $1,\ldots ,n$. Due to  the equality of the right hand sides in \eqref{re1} and in \eqref{re2}, we can use crossing symmetry property \eqref{csym_equiv} to express \eqref{xx1} as
\begin{align}
&\tr_{1,\ldots ,n} \ts A^{(n)}   \ts A\cdotrl \left( \left( C\ts  \Tc_{[n]}(u_{[n]})  \ts D\ts  \Tc_0 (u_0)\right)\cdot B\right)
,\quad\text{where}\quad A=\wvr{R}_{1n}^{12}(u_0|u_{[n]}),\label{xx4}\\
& B=\wvr{R}_{1n}^{12}(u_0|u_{[n]})^{-1},\quad C=\wvr{R}_{1n}^{12\ts r}(u_0|u_{[n]})^{-1}\fand D= \wvr{R}_{1n}^{12\ts r}(u_0+hN |u_{[n]}).\non
\end{align}
Recall that the meaning of the superscript r is explained in the statement of Lemma \ref{superrr}.  
By \eqref{xxx1}  and \eqref{xxx2} we have
\beq\label{xx3}
 A^{(n)}  \ts Z\ts=\ts \cev{Z}\ts  A^{(n)}\quad\text{for }Z=A ,B ,C ,D .
\eeq
Therefore, since $\left(A^{(n)}\right)^2=A^{(n)}$, we conclude that \eqref{xx4} is equal to
$$
\tr_{1,\ldots ,n} \ts    \cev{A}\cdotrl \left( \left( \cev{C}\ts \left(A^{(n)}\right)^2 \ts \Tc_{[n]}(u_{[n]})  \ts D\ts  \Tc_0 (u_0)\right)\cdot B\right).
$$
Next, we employ \eqref{xxx3} and \eqref{xx3}  to move one copy of $A^{(n)}$ to the left and another copy of $A^{(n)}$ to the right, thus getting
$$
\tr_{1,\ldots ,n}     \ts A^{(n)}\ts A \cdotrl \left( \left( C  \ts  \cev{\Tc}_{[n]}(u_{[n]})  \ts \cev{D} \ts  \Tc_0 (u_0)\right)\cdot \cev{B}\ts    A^{(n)} \right) .
$$
By the cyclic property of the trace and $\left(A^{(n)}\right)^2=A^{(n)}$ this equals to
$$
\tr_{1,\ldots ,n}     \ts  A \cdotrl \left( \left( C  \ts  \cev{\Tc}_{[n]}(u_{[n]})  \ts \cev{D} \ts  \Tc_0 (u_0)\right)\cdot \cev{B}\ts    A^{(n)} \right) .
$$
Finally,  using \eqref{xxx3} and \eqref{xx3} we move the remaining copy of $A^{(n)}$ to the left:
$$
\tr_{1,\ldots ,n}     \ts   A \cdotrl \left( \left( C  \ts  A^{(n)}\ts   \Tc_{[n]}(u_{[n]})  \ts D \ts  \Tc_0 (u_0)\right)\cdot B  \right).
$$

Since $B \ts A =1$,  by employing the cyclic property of the trace and moving the tensor factors $1,\ldots ,n$ of $A$ to the right we obtain
$$
\tr_{1,\ldots ,n}      \left( C  \ts A^{(n)}\ts  \Tc_{[n]}(u_{[n]})  \ts D \ts  \Tc_0 (u_0)\right)\cdot\left(B\ts A\right)\ts
=\ts
\tr_{1,\ldots ,n}   \ts   C  \ts A^{(n)}\ts  \Tc_{[n]}(u_{[n]})  \ts D \ts  \Tc_0 (u_0).
$$
In order to finish the proof, i.e. to verify   \eqref{center2}, it is sufficient to check that 
\beq\label{xx7}
\tr_{1,\ldots ,n} \ts      C \ts  A^{(n)} \ts  \Tc_{[n]}(u_{[n]})  \ts D 
\eeq
  is equal to   
$\Tf_{[n]} (u)$. However,  crossing symmetry
property \eqref{csym_equiv} implies $C  \cdotlr D  =1$.
Hence, using the cyclic property of the trace and moving the tensor factors $1,\ldots ,n$ of $C$ to the right, we rewrite \eqref{xx7} as follows:
\begin{align*}
&\tr_{1,\ldots ,n} \ts C \ts  A^{(n)} \ts  \Tc_{[n]}(u_{[n]})  \ts D =
\tr_{1,\ldots ,n} \ts C \cdotlr \left(A^{(n)}  \Tc_{[n]}(u_{[n]})  \ts D\right)\\
=&
\tr_{1,\ldots ,n}    \left( A^{(n)} \ts \Tc_{[n]}(u_{[n]})    \right)\cdot ( C  \cdotlr D  )
= \tr_{1,\ldots ,n} \ts  A^{(n)} \ts\Tc_{[n]}(u_{[n]})    = \Tf_{[n]}(u),
\end{align*}
thus proving the theorem.
\end{prf}

\begin{rem}\label{detremark}
Recall   action \eqref{homohomo} of the   quantum current algebra on the vacuum module.
It is worth noting that the action of $\Tf_{[n]} (u)$  coincides with  the action of    certain Laurent  series $\wtld{\TT}_{1^n}(u)$, whose coefficients belong to the center of the completed double Yangian $\wtld{\DY}_{-N}(\gl_N)$ 
at the critical level $c=-N$; see \cite[Theorem 4.4]{JKMY} for more details. 
\end{rem}

\begin{rem}
In \cite[Proposition 3.1]{D}, Ding found a realization of the quantum affine algebra $U_q (\widehat{\gl}_N)$ via commutation relation of the form \eqref{re}, which  involves trigonometric $R$-matrix. The proof of Theorem \ref{centralthm} and  fusion procedure from \cite{C} suggest that the analogous construction of central elements of the completed quantum affine algebra $U_q (\widehat{\gl}_N)$ at the critical level  might be  given in terms of Ding's realization, as long as   the corresponding (trigonometric) $R$-matrix is appropriately normalized.
The image of such family of central elements in the $RTT$ presentation of the completed quantum affine algebra $U_q (\widehat{\gl}_N)$ at the critical level  should coincide with the coefficients of the Laurent series $\ell_k(z)$, $k=1,\ldots ,N$, as defined in \cite[Section 3]{FJMR}, thus providing a new proof of \cite[Theorem 3.2]{FJMR}.
\end{rem}

\section{Vacuum module over the  quantum current algebra } \label{sec3}
In this section, we introduce the vacuum module  $\Vc_c (\R)$ for the quantum current algebra  and we show that it is isomorphic, as a $\CC[[h]]$-module, to the vacuum module over the double Yangian. This allows us to employ Etingof--Kazhdan's construction and obtain the structure of quantum vertex algebra on $\Vc_c (\R)$. Next, in parallel with representation theory of the affine  Lie algebras, we introduce the notion of restricted module for the  algebra $\A(\R)$. Finally, we prove that restricted level $c$ modules for the   algebra $\A(\R)$ are  naturally equipped with a
structure of module  for the quantum vertex algebra $\Vc_c (\R)$ and vice versa.

\subsection{Properties of the vacuum module}\label{sec31} 

Introduce the subset $W_2^c$ of  $A_2^c$ by
$$W_2^c = \left\{\tz_{i_1\ts j_1\ts\ldots\ts i_n\ts j_n }^{(r_1,\ldots ,r_n)} \in A_2^c\,:\,r_k\geqslant 0\text{ for some }k=1,\ldots,n  \right\}.$$
 Let $\W_c (\R)$ be the left ideal in the  algebra $\A_c (\R)$ generated by the set $W_2^c$. 
Introduce the completion of $\W_c (\R)$ as the inverse limit
$$\W_c' (\R) =\lim_{\longleftarrow} \ts \W_c (\R)\ts /\ts \W_c (\R)\cap \I_p (\R).$$
Then the $h$-adic completion $\wtld{[\W_c' (\R)]}$ of
$$[\W_c' (\R)]=\left\{a\in \A_c (\R)\,:\, h^n a\in \W_c' (\R)\text{ for some }n\geqslant 0\right\}$$
is also a left ideal in $\A_c (\R)$.
Define the {\em vacuum module} $\V_c(\R)$  as the   quotient of the algebra $\A_c (\R)$ by its left ideal  $\wtld{[\W_c' (\R)]}$,   
\beq\label{quotient2}
 \V_c (\R)\,=\,\A_c (\R)\,/\, \wtld{[\W_c' (\R)]}.
\eeq
Observe that the canonical map $\A (\R)\to \V_c (\R)$ maps  the left ideal $\I_p^h(\R)$   to $h^p  \V_c (\R)$.

\begin{pro}\label{free2}
The vacuum module $\V_c (\R)$ is topologically free.
\end{pro}

\begin{prf}
The algebra $\A_c (\R)$ is topologically free, so the proposition can be verified by arguing as in the proof of Proposition \ref{free}.
\end{prf}

Denote the image of the unit $1\in\A_c (\R)$ in  quotient \eqref{quotient2} by $\vac$.
Let $V^c$ be the set of   all elements $a\cdot \vac \in \V_c(\R)$ such that $a\in A_2^c\setminus W_2^c$, i.e.
$$V^c=\left\{\vac\right\}\cup\left\{\tz_{i_1\ts j_1\ts\ldots\ts i_n\ts j_n }^{(r_1,\ldots ,r_n)}\vac  \,:\,\tz_{i_1\ts j_1\ts \ldots\ts i_n\ts j_n }^{(r_1,\ldots ,r_n)} \in A_2^c \text{ and }r_k < 0\text{ for all }k=1,\ldots,n \right\} .$$
As a direct consequence of Proposition \ref{pro_1} we obtain
\begin{pro}\label{pro_2}
The $\CC[[h]]$-module $\V_c(\R)$ coincides with the $h$-adic completion of
$$[\xpan V^c]=\left\{v\in \V_c(\R)\,:\, h^n v\in \xpan V^c\text{ for some }n\geqslant 0\right\}.$$
\end{pro}

We are now ready to prove the main result in this subsection.

\begin{thm}\label{izomorfizam}
The assignments 
\beq\label{mhomo}
\Tc_{[n]}(u)\vac\,\mapsto\, T_{[n]}^+ (u)\vac
\eeq
with $n\geqslant 1$ and the variables $u=(u_1,\ldots ,u_n)$ define a   $\CC[[h]]$-module isomorphism
\beq\label{mhomo2}
\V_c(\R)\,\to\, \Vc_{c}(\gl_N).
\eeq
\end{thm}

\begin{prf}
In order to prove that \eqref{mhomo} defines a  homomorphism of $\CC[[h]]$-modules,
it is sufficient to check that the elements of the set $W_2^c$ belong to the kernel of   $\CC[[h]]$-module map \eqref{chhomo}. Let
$\tz=\tz_{i_1\ts j_1\ts\ldots\ts i_n\ts j_n }^{(r_1,\ldots ,r_n)}$ be an arbitrary element of $W_2^c$. Then $r_k\geqslant 0$ for some $k=1,\ldots ,n$. The image $\wht{\tz}\in \Vc_{c}(\gl_N)$ of the element $\tz$, with respect to map \eqref{chhomo}, coincides with the coefficient of the variables $u_1^{-r_1 -1}\ldots u_n^{-r_n -1}$ in the matrix entry $e_{i_1\ts j_1}\ot\ldots\ot e_{i_n\ts j_n}$ of the expression 
\beq\label{expr_675}
T_{[n]}^+ (u)T_{[n]} (u+hc/2)^{-1}\vac=T_{[n]}^+ (u)\vac\,\in\,(\ndo\CC^N)^{\ot n}\ot \Vc_{c}(\gl_N)[[u_1,\ldots ,u_n]]. 
\eeq
Since \eqref{expr_675} does not contain any negative powers of the variable $u_k$, we conclude that $\wht{\tz}$ equals $0$, as required.
Therefore, assignments \eqref{mhomo} define a $\CC[[h]]$-module homomorphism
$\V_c(\R)\to \Vc_{c}(\gl_N)$.
Moreover, by the Poincar\'e--Birkhoff--Witt theorem for the double Yangian, see \cite[Theorem 2.2]{JKMY}, this map   is surjective.

Let us prove that the assignments 
\beq\label{inverz}
T_{[n]}^+ (u)\vac\,\mapsto\,\Tc_{[n]}(u)\vac
\eeq
with $n\geqslant 1$ and the variables $u=(u_1,\ldots ,u_n)$ define
  a   $\CC[[h]]$-module epimorphism
\beq\label{inverz2}
\Vc_{c}(\gl_N)\,\to\, \V_c(\R).
\eeq
In order to verify that  \eqref{inverz}  defines a  $\CC[[h]]$-module homomorphism, it is sufficient to check that \eqref{inverz} maps the ideal of dual Yangian relations \eqref{RTT1}  to itself. However, this follows from Lemma \ref{newl}. Indeed, by \eqref{RTT1} we have
\beq\label{new3}
\R_{i\ts i+1}(u_i -u_{i+1}) T_{[n]}^+(u )\vac \R_{i\ts i+1}(u_i -u_{i+1})^{-1}
=    P_{i\ts i+1}T_{[n]}^+(u_{i\leftrightarrow i+1} )\vac    P_{i\ts i+1}
\eeq
 for any $i=1,\ldots, n-1$.
The images of the left and the right hand side in \eqref{new3}, with respect to \eqref{inverz}, are equal to the left and the right hand side in the equality
$$\R_{i\ts i+1}(u_i -u_{i+1}) \Tc_{[n]}(u)\vac \R_{i\ts i+1}(u_i -u_{i+1})^{-1}
=P_{i\ts i+1} \Tc_{[n]}(u_{i\leftrightarrow i+1})\vac    P_{i\ts i+1},
$$
which  follows by applying \eqref{new} on $\vac\in\V_c(\R)$. Hence, we conclude that \eqref{inverz} defines a $\CC[[h]]$-module homomorphism.
Moreover, the set $V^c$ is contained within the image of  \eqref{inverz2}, so  the  map   is surjective by Proposition \ref{pro_2}.
Finally, since maps \eqref{mhomo2} and \eqref{inverz2} are inverses of each other, the theorem follows.
\end{prf}

We now proceed towards the definition of restricted $\A (\R)$-module, which is motivated by the following proposition.

\begin{pro}\label{baka}
For any integer $n\geqslant 1$ and the variables $u=(u_1,\ldots ,u_n)$ we have
$$\Tc_{[n]}(u)w\,\in\,(\ndo\CC^N)^{\ot n}\ot \V_c(\R)((u_1,\ldots,u_n))[[h]]\quad\text{for all } w\in\V_c(\R).$$
\end{pro}

\begin{prf}
By applying \eqref{re1} on the element $\vac\in\Vc_c(\R)$ we get
$$\Tc_{[n+m]}(u,v)\vac=\Tc_{[n]}^{13}(u) \,\R_{nm}^{12}(u+hc|v)^{-1} \,\Tc_{[m]}^{23}(v)\,\R_{nm}^{12}(u|v)\vac$$
for the variables $u=(u_1,\ldots ,u_n)$ and $v=(v_1,\ldots ,v_m)$.
Note that this expression contains only nonnegative powers of the variables $u_1,\ldots ,u_n$ and $v_1,\ldots ,v_m$. Using crossing symmetry property \eqref{csym_equiv} we   move all $R$-matrices to the left hand side, thus getting
$$\R_{nm}^{12}(u+hc+hN|v) \cdotrl\left(\Tc_{[n+m]}(u,v)\vac\cdot\R_{nm}^{12}(u|v)^{-1}\right)=\Tc_{[n]}^{13}(u) \left(\Tc_{[m]}^{23}(v)\vac\right).$$
Finally, we observe that, for any given integers $a_1,\ldots ,a_m\geqslant 0$ and $p\geqslant 1$, the coefficient of the monomial $v_1^{a_1}\ldots v_m^{a_m}$ on left hand side of the given equation contains only finitely many negative powers of the variables $u_1,\ldots ,u_n$ modulo $h^p$. Since the set of coefficients of the matrix entries of all $\Tc_{[m]}(v)\vac$ with $m\geqslant 0$ coincides with $V^c$, the proposition   follows from Proposition \ref{pro_2}. 
\end{prf}

An  $\A(\R)$-module $W$ is said to be {\em restricted} if
 $W$ is a topologically free $\CC[[h]]$-module such that   $$\Tc (u)w\,\in\,\ndo\CC^N\ot W ((u))[[h]]\quad\text{for all }w\in W.$$
Also, as usual, an
$\A(\R)$-module $W$ is said to be a {\em level $c$} module if the central element $C\in \A(\R)$ acts on $W$ as a  scalar multiplication by some $c\in\CC$. 
Propositions \ref{free2} and  \ref{baka} imply that the vacuum module $\V_c(\R)$ is  restricted $\A(\R)$-module of level $c$.
 
\begin{pro}\label{restricted_important_very}
Let $W$ be a restricted $\A(\R)$-module. Then
\beq\label{restricted_important}
\Tc_{[n]}(u)w \,\in\, (\ndo\CC^N)^{\ot n}\ot W((u_1,\ldots,u_n))[[h]]\quad\text{for all }w\in W\text{ and }n\geqslant 1.
\eeq
\end{pro}

\begin{prf}
The statement follows by  induction on $n$ which is based on   relations  \eqref{re1}--\eqref{re2} and arguments from the proof of Proposition \ref{teemlema}. 
\end{prf}

Suppose $W$  is a  restricted $\A(\R)$-module of level $c$. For the variable $z$ and the variables $u=(u_1,\ldots ,u_n)$   introduce the  elements of 
$(\ndo\CC^N)^{\ot n}\ot (\ndo W) ((z))[[u_1,\ldots,u_n,h]]$ by
$$\Tc_{[n]}(u|z)=\hspace{-4pt}\prod_{i=1,\ldots ,n}^{\longrightarrow}\hspace{-4pt} \left(\Tc_{i}(z+u_i)\R_{i\ts i+1}(u_i -u_{i+1}+hc)^{-1}\ldots \R_{i\ts n}(u_i -u_{n}+hc)^{-1}  \right)\,\cdot\, \cev{\R}_{[n]}(u).
$$
In particular, we have $\Tc_{[1]}(u|z)=\Tc(z+u)$.
The following proposition  is required in the proof of Theorem \ref{main}.

\begin{pro}\label{reflectionprop2}
Let $W$ be a restricted $\A(\R)$-module of level $c$.
\begin{enumerate}[(a)]
\item\label{stat1}
For any integers $n,m\geqslant 1$   the  equalities
\begin{align}
&\Tc_{[n+m]}(z_1+u_1,\ldots , z_1+u_n,z_2+v_1,\ldots ,z_2 +v_m)\non\\
=&\Tc_{[n]}^{13}(u|z_1 ) \,\R_{nm}^{12}(u|v|z_1-z_2+hc)^{-1} \,\Tc_{[m]}^{23}(v|z_2)\,\R_{nm}^{12}(u|v|z_1 -z_2)\label{re_z}
\\
=&\R_{nm}^{12 }(u|v|-z_2+z_1)^{-1}\,\Tc_{[m]}^{23}(v|z_2) \,\R_{nm}^{12 }(u|v|-z_2 +z_1-hc)\,\Tc_{[n]}^{13}(u|z_1)\non
\end{align}
hold
in
$$
(\ndo\CC^N)^{\ot n}\ot (\ndo\CC^N)^{\ot m}\ot (\ndo W) ((z_1,z_2))[[u_1,\ldots,u_n,v_1,\ldots, v_m, h]].$$
\item\label{stat2} 
For any $n\geqslant 2$ and $i=1,\ldots ,n-1$ the equality
\beq\label{new_z}
\R_{i\ts i+1}(u_i -u_{i+1}) \Tc_{[n]}(u|z) \R_{i\ts i+1}(u_i -u_{i+1})^{-1}
=P_{i\ts i+1} \Tc_{[n]}(u_{i\leftrightarrow i+1}|z)    P_{i\ts i+1}
\eeq
holds  in
$(\ndo\CC^N)^{\ot n}\ot (\ndo W) ((z,u_{i}))[[u_1,\ldots,u_n, h]].$
\end{enumerate}
\end{pro}

\begin{prf}
Proposition \ref{reflectionprop} implies the first and
Lemma \ref{newl} implies the second statement of the proposition. 
\end{prf}

\subsection{Vacuum module as a quantum vertex algebra}\label{sec33}

From now on, we  often identify the $\CC[[h]]$-modules $\V_c(\R)$ and $\Vc_{c}(\gl_N)$ via the $\CC[[h]]$-module isomorphism established in Theorem \ref{izomorfizam}. For example, we utilize such identification in the next theorem, which is due to Etingof and Kazhdan; see  \cite[Theorem 2.3]{EK}. However, the vertex operator map in \eqref{qva1} is expressed somewhat differently from the original version in \cite{EK}, so we demonstrate in the proof that both definitions coincide. 
\begin{thm}\label{EK:qva}
For any $c\in \CC$
there exists a unique  structure of quantum vertex algebra
on  $\V_c(\R)$ such that the vacuum vector is
$\vac\in \V_c(\R)$, the vertex operator map is defined by
\beq\label{qva1}
Y\big(T_{[n]}^+ (u)\vac,z\big)=\Tc_{[n]}(u|z), 
\eeq
  the map $D$ is defined by
\beq\label{mapdi}
e^{zD} T_{[n]}^+ (u)\vac = T_{[n]}^+ (u|z)\vac 
\eeq
  and the map $\mathcal{S}(z)$ is defined by  
\begin{align}
\mathcal{S}(z)&\Big(\overline{R}_{nm}^{  12}(u|v|z)^{-1}  T_{[m]}^{+24}(v) 
\overline{R}_{nm}^{  12}(u|v|z-h  c)  T_{[n]}^{+13}(u)(\vac\otimes \vac) \Big)\nonumber\\
 =&T_{[n]}^{+13}(u)  \overline{R}_{nm}^{  12}(u|v|z+h  c)^{-1} 
 T_{[m]}^{+24}(v)  \overline{R}_{nm}^{  12}(u|v|z)(\vac\otimes \vac)\label{qva2}
\end{align}
for operators on
$
(\ndo\mathbb{C}^{N})^{\otimes n} \otimes
(\ndo\mathbb{C}^{N})^{\otimes m}\otimes \V_c(\R) \ot \V_c(\R) $.
\end{thm}

\begin{prf}
Since maps \eqref{mapdi} and \eqref{qva2} coincide with the original maps from \cite[Theorem 2.3]{EK},   we only have to check that \eqref{qva1} coincides with the original definition of the vertex operator map in \cite{EK}, which is given by
\beq\label{ekqva1}
 T_{[n]}^+ (u)\vac\,\mapsto\, T^+_{[n]}(u|z)T_{[n]}(u|z+hc/2)^{-1}\quad\text{for }n\geqslant 1. 
\eeq
Due to  Poincar\'e--Birkhoff--Witt theorem for the double Yangian  \cite[Theorem 2.2]{JKMY},
it is sufficient to prove that for any integer $m\geqslant 1$ the actions of \eqref{qva1} and \eqref{ekqva1} on $\Tc_{[m]}(v)\vac\equiv T_{[m]}^+ (v)\vac$, where $v=(v_1,\ldots ,v_m)$, coincide. 
By applying \eqref{ekqva1} on $T_{[m]}^+ (v)\vac$ and using      relation \eqref{rtt3}, together with crossing symmetry  property \eqref{csym_equiv},    we obtain
\beq\label{ekqva2}
T^{+13}_{[n]}(u|z)T^{13}_{[n]}(u|z+hc/2)^{-1} T_{[m]}^{+23} (v)\vac= A\cdotrl \left(T_{[n]}^{+13} (u|z)T_{[m]}^{+23} (v)\vac\cdot B\right),
\eeq
where
$$A=\R_{nm}^{12}(u|v|z+hc+hN)\fand B=\R_{nm}^{12}(u|v|z)^{-1}.$$
On the other hand, by applying \eqref{qva1} on $\Tc_{[m]}(v)\vac$ and using      commutation relation \eqref{re_z}, together with crossing symmetry property \eqref{csym_equiv},  we get
\beq\label{xqva}
\Tc_{[n]}^{13}(u|z)\Tc_{[m]}^{23}(v)\vac=A\cdotrl \left(\Tc_{[n+m]}(z+u_1,\ldots,z+u_n,v_1,\ldots ,v_m)\vac\cdot B\right).
\eeq
Since  $\CC[[h]]$-module isomorphism \eqref{mhomo2} maps 
the right hand side of  \eqref{xqva} to the right hand side of \eqref{ekqva2},
definitions in \eqref{qva1} and in \eqref{ekqva1} coincide, as required.
\end{prf}

\subsection{Main result}\label{subsec34}
In the following lemma, we introduce certain map $\mc (z)$ which is used in the proof of Theorem \ref{main}.

\begin{lem}\label{lemmam}
For any $c\in\CC$ the assignments
\beq\label{m2m2m2}
T_{[n]}^{+13 }(u)T_{[m]}^{+24 }(v)(\vac\ot\vac)
\, \mapsto\,
\R_{nm}^{12}(u|v|z)^{-1}T_{[m]}^{+23 }(v)
\R_{nm}^{12}(u|v|z-hc)
T_{[n]}^{+14 }(u) (\vac\ot\vac)
\eeq
with $n,m\geqslant 1$ and the variables $u=(u_1,\ldots ,u_n) $ and  $v=(v_1,\ldots ,v_m) $ define a
 $\CC[[h]]$-module map 
\beq\label{m2}
\mc (z)\,\colon\,\Vc_c(\R)\ot\Vc_c(\R)\,\to\, \Vc_c(\R)\ot\Vc_c(\R)\ot\CC((z)).
\eeq
Moreover, the following equalities hold:  
\begin{align}
 &Y(z)\left( \mc(-z) \left( T_{[n]}^{+13 }(u)T_{[m]}^{+24 }(v)(\vac\ot\vac)\right)\right) =T_{[n]}^{+13 }(u)T_{[m]}^{+23 }(v|z)\vac ;\label{lemmaa}\\
& \Sc(z)\left(P'\left( \mc(z)\left( T_{[n]}^{+13 }(u)T_{[m]}^{+24 }(v)(\vac\ot\vac)\right)\right)\right) \non\\
&\qquad\qquad=T_{[n]}^{+13}(u)  \overline{R}_{nm}^{  12}(u|v|z+h  c)^{-1} 
 T_{[m]}^{+24}(v)  \overline{R}_{nm}^{  12}(u|v|z)(\vac\ot\vac)\label{lemmab}
\end{align}
for operators on
$
(\ndo\mathbb{C}^{N})^{\otimes n} \otimes
(\ndo\mathbb{C}^{N})^{\otimes m}\otimes \V_c(\R) \ot \V_c(\R), $
where $P'\colon w_1\ot w_2\mapsto w_2\ot w_1$ denotes the permutation operator on $\V_c(\R) \ot\V_c(\R)$. 
\end{lem}

\begin{prf}
The fact that $\CC[[h]]$-module map \eqref{m2} is well-defined can be proved by a simple calculation which relies on Yang--Baxter equation \eqref{ybe} and defining relations \eqref{RTT1} for the dual Yangian. The proof of equalities \eqref{lemmaa} and \eqref{lemmab} is also straightforward. The former employs  unitarity property \eqref{uni} and relations \eqref{rtt1} and \eqref{rtt3} while the latter follows directly from \eqref{qva2}. 
\end{prf}

The following theorem is the main result of this paper.
\begin{thm}\label{main}
Let $W$ be a restricted $\A(\R)$-module of level $c\in\CC$. There exists a unique structure of $\V_c(\R)$-module on $W$ satisfying
\beq\label{main1}
Y_W(T_{[n]}^+ (u)\vac,z)=\Tc_{[n]}(u|z)
\eeq
for all $n\geqslant 1$.
Conversely, let $(W,Y_W)$ be a $\V_c(\R)$-module for some $c\in\CC$. There exists a unique structure of restricted $\A(\R)$-module of level $c$ on $W$ satisfying
\beq\label{main2}
\Tc(z)=Y_W(T^+(0)\vac,z).
\eeq
Moreover, a topologically free $\CC[[h]]$-submodule $W_1$ of $W$ is a $\V_c(\R)$-submodule of $W$ if and only if  $W_1$ is an $\A(\R)$-submodule of $W$.
\end{thm}

\begin{prf}
Let $W$ be a restricted $\A(\R)$-module of level $c\in\CC$. Equalities \eqref{main1} with $n\geqslant 1$, together with $Y_W(\vac,z)=1_W
$, define a $\CC[[h]]$-module map $Y_W(z)\colon V\ot W\to W((z))[[h]]$. Indeed, the fact that $Y_W(z)$ is well-defined can be verified by a simple calculation which relies on defining relations \eqref{RTT1} for the dual Yangian and   \eqref{new_z}.

We  now prove that the map $Y_W(z)$, as defined by \eqref{main1}, satisfies weak associativity \eqref{associativityw}. 
Let  $w$ be an arbitrary element of $W$.  Consider the expression
\beq\label{xtemp2}
Y_W(T_{[n]}^{+\ts 13} (u)\vac,z_0 +z_2) Y_W(T_{[m]}^{+\ts 23} (v)\vac,z_2)w ,
\eeq
which corresponds to the first summand in \eqref{associativityw}.
By \eqref{main1}, the given expression   equals
\beq\label{xtemp2222}
\Tc_{[n]}^{13}(u|z_0 +z_2)\Tc_{[m]}^{23}(v|z_2)w.
\eeq
By combining  \eqref{re_z}   and crossing symmetry property \eqref{csym_equiv} we   express \eqref{xtemp2222} as
\beq\label{xtemp3}
\R_{nm}^{12}(u|v|z_0+hc+hN)\cdotrl\left(\Tc_{[n+m]}(x)w\cdot\R_{nm}^{12}(u|v|z_0)^{-1}\right),
\eeq
where the variables $x=(x_1,\ldots ,x_{n+m})$ are given by
$$x=(z_0+z_2+u_1,\ldots ,z_0+z_2+u_n,z_2+v_1,\ldots ,z_2+v_m).$$

Next, consider the expression
\beq\label{xtemp4}
Y_W(Y(T_{[n]}^{+\ts 13} (u)\vac,z_0)T_{[m]}^{+\ts 23} (v)\vac,z_2)w,
\eeq
which corresponds to the second summand in \eqref{associativityw}.
By using \eqref{qva1} and the identification $T_{[m]}^{+} (v)\vac\equiv \Tc_{[m]}^{+} (v)\vac$  we express \eqref{xtemp4} as
\beq\label{xtemp5}
Y_W\left(\Tc_{[n]}^{13}(u|z_0 )\Tc_{[m]}^{23}(v)\vac,z_2\right)w.
\eeq
As before, we employ   \eqref{re_z} and crossing symmetry property \eqref{csym_equiv} to write \eqref{xtemp5} as
\beq\label{xtemp6}
Y_W\left(\R_{nm}^{12}(u|v|z_0+hc+hN)\cdotrl (\Tc_{[n+m]}(y')\vac\cdot\R_{nm}^{12}(u|v|z_0)^{-1} ),z_2\right)w,
\eeq
where the variables $y'=(y'_1,\ldots,y'_{n+m})$ are given by
$$y'=(z_0+u_1,\ldots ,z_0+u_n,v_1,\ldots ,v_m).$$
Due to \eqref{main1}, the expression in \eqref{xtemp6} is equal to
\beq\label{xtemp7}
\R_{nm}^{12}(u|v|z_0+hc+hN)\cdotrl\left(\Tc_{[n+m]}(y)w\cdot\R_{nm}^{12}(u|v|z_0)^{-1}\right),
\eeq
where the variables $y=(y_1,\ldots ,y_{n+m})$ are given by
$$y=(z_2+z_0+u_1,\ldots ,z_2+z_0+u_n,z_2+v_1,\ldots ,z_2+v_m).$$

Observe that   \eqref{xtemp3}  and \eqref{xtemp7} are not equal. Indeed, due to our expansion convention from Section \ref{sec13}, the former is expanded in nonnegative powers of the variable $z_2$ while the latter is expanded in nonnegative powers of the variable $z_0$.
Fix arbitrary nonnegative integers $k,r_1,\ldots,r_{n},s_1,\ldots ,s_m$. 
Since $W$ is a  restricted $\A(\R)$-module, Proposition \ref{restricted_important_very} implies that the coefficient of the monomial $u_1^{r_1}\ldots u_n^{r_n}v_1^{s_1}\ldots v_m^{s_m}$ in 
\beq\label{xtemp8}
\Tc_{[n+m]}(z+u_1,\ldots ,z+u_n,v_1,\ldots v_m)w
\eeq
possesses only finitely many negative powers of the variable $z$ modulo $h^k$. Choose an integer $r\geqslant 0$ such that the coefficient of $u_1^{r_1}\ldots u_n^{r_n}v_1^{s_1}\ldots v_m^{s_m}$ in 
$$z^r\ts \Tc_{[n+m]}(z+u_1,\ldots ,z+u_n,v_1,\ldots v_m)w$$
 possesses only nonnegative powers of the variable $z$ modulo $h^k$.\footnote{Notice that the integer $r\geqslant 0$ also depends on the choice of $w\in W$.} 
Then, by our discussion, the coefficients
of the monomial $u_1^{r_1}\ldots u_n^{r_n}v_1^{s_1}\ldots v_m^{s_m}$ in
$$(z_0 +z_2)^r\ts Y_W(T_{[n]}^{+\ts 13} (u)\vac,z_0 +z_2) Y_W(T_{[m]}^{+\ts 23} (v)\vac,z_2)w $$
and in
$$(z_0 +z_2)^r\ts Y_W(Y(T_{[n]}^{+\ts 13} (u)\vac,z_0)T_{[m]}^{+\ts 23} (v)\vac,z_2)w$$
 coincide modulo $h^k$, so we conclude that weak associativity \eqref{associativityw} holds.
Hence, Lemma \ref{usefullemma}  implies that $(W,Y_W)$ is a $\Vc_c(\R)$-module.

Conversely, 
let $(W,Y_W)$ be a $\V_c (\R)$-module. For an arbitrary $w\in W$ let us apply Jacobi identity \eqref{wjacobi} on the last three tensor factors of the expression
\begin{align}
A(z_0)&\coloneqq \mc(-z_0)(T^{+}_{13}(u)  T^+_{24}(v)(\vac\ot\vac))\ot w\label{azenula2}\\
&=\R (-z_0+u-v)^{-1}T_{23}^{+ }(v)
\R (-z_0+u-v-hc)
T_{14}^{+ }(u) (\vac\ot\vac) \ot w\label{azenula},
\end{align}
which belongs to $(\ndo\CC^N)^{\ot 2}\ot \V_c(\R)^{\ot 2}\ot W ((z_0))[[u,v,h]]$. 
By applying the first term in \eqref{wjacobi}, 
$$z_0^{-1}\delta\left(\frac{z_1 -z_2}{z_0}\right) Y_W(z_1)(1\ot Y_W(z_2))$$
on  $A(z_0)$, as given in  \eqref{azenula}, we get
\begin{align}
&z_0^{-1}\delta\left(\frac{z_1 -z_2}{z_0}\right)\R (-z_1+z_2+u-v)^{-1}Y_W(T^{+ }_2(v)\vac,z_1)\non\\
&\qquad\times\R (-z_1+z_2+u-v-hc)
Y_W(T^{+ }_1(u)\vac,z_2) w.\label{m1}
\end{align}
By applying the second term in \eqref{wjacobi}, 
$$-z_0^{-1}\delta\left(\frac{z_2-z_1}{-z_0}\right) Y_W(z_2)(1\ot Y_W(z_1)) \left(\Sc(-z_0)P'\ot 1\right)$$
on $A(z_0)$, as given in  \eqref{azenula2}, and using  \eqref{lemmab}  we get
\begin{align}
&-z_0^{-1}\delta\left(\frac{z_2-z_1}{-z_0}\right) Y_W(T^{+ }_1(u)\vac,z_2) \non\\
&\qquad\times
\R (z_2-z_1+u-v+hc)^{-1}
Y_W(T^{+ }_2(v)\vac,z_1) \R (z_2-z_1+u-v)   w.
\label{m3}
\end{align}
Finally, by applying the third term in \eqref{wjacobi}, 
$$z_2^{-1}\delta\left(\frac{z_1 -z_0}{z_2}\right)Y_W(z_2)(Y(z_0)\ot 1 )$$
on $A(z_0)$, as given in  \eqref{azenula2}, and using   \eqref{lemmaa}  we get
\begin{align}
z_2^{-1}\delta\left(\frac{z_1 -z_0}{z_2}\right)Y_W(T_{1}^+(u)T^{+}_{2}(z_0+v)\vac,z_2 )w.\label{m4}
\end{align}

Since \eqref{m4} does not contain any negative powers of the variable $z_0$, the sum of residues of expressions \eqref{m1} and \eqref{m3}, with respect to the variable $z_0$, equals $0$. Therefore, by taking the residue $\rez_{z_0}$ we obtain
\begin{align*}
&\R (-z_1+z_2+u-v)^{-1}Y_W(T^{+ }_2(v)\vac,z_1)\R (-z_1+z_2+u-v-hc)
Y_W(T^{+ }_1(u)\vac,z_2) w\\
=&
Y_W(T^{+ }_1(u)\vac,z_2) 
\R (z_2-z_1+u-v+hc)^{-1}
Y_W(T^{+ }_2(v)\vac,z_1) \R (z_2-z_1+u-v)   w.
\end{align*}
Both sides of this  equality   contain only nonnegative powers of the variables $u$ and $v$. In particular, their constant terms, with respect to $u$ and $v$, coincide, i.e. we have
\begin{align*}
&\R (-z_1+z_2 )^{-1}Y_W(T^{+ }_2(0)\vac,z_1)\R (-z_1+z_2 -hc)
Y_W(T^{+ }_1(0)\vac,z_2) w\\
=&
Y_W(T^{+ }_1(0)\vac,z_2) 
\R (z_2-z_1+hc)^{-1}
Y_W(T^{+ }_2(0)\vac,z_1) \R (z_2-z_1)   w.
\end{align*}
Since the element $w\in W$ was arbitrary, the following equality holds on $W$:
\begin{align}
&\R (-z_1+z_2 )^{-1}Y_W(T^{+ }_2(0)\vac,z_1)\R (-z_1+z_2 -hc)
Y_W(T^{+ }_1(0)\vac,z_2)  \non\\
=&
Y_W(T^{+ }_1(0)\vac,z_2) 
\R (z_2-z_1+hc)^{-1}
Y_W(T^{+ }_2(0)\vac,z_1) \R (z_2-z_1)    \label{reW}.
\end{align}
Observe that \eqref{reW} coincides with commutation relation \eqref{re} at $C=c$.
Therefore, since  $Y_W(T^+(0)\vac,z)w$ belongs to $\ndo\CC^N \ot W ((z))[[h]]$ for all $w\in W$, we conclude that \eqref{main2} defines a structure of a level $c$ restricted $\A(\R)$-module on $W$.

Finally, let us prove the last assertion of the theorem. Suppose that $W_1$ is a $\V_c(\R)$-submodule of $W$. Then for any $w_1\in W_1$ we have
$$\Tc (z)w_1=Y_W(T^+(0)\vac,z)w_1 \,\in\, \ndo\CC^N \ot W_1 ((z))[[h]],$$
so $W_1$ is an $\A(\R)$-submodule of $W$. Conversely, suppose that $W_1$ is a topologically free $\A (\R)$-submodule of $W$.
Clearly, $W_1$ is a restricted $\A (\R)$-module (of level $c$), so Proposition \ref{restricted_important_very} implies that for any $w_1\in W_1$ we have
$$
\Tc_{[n]}(v)w_1\,\in\,\left(\ndo\CC^N\right)^{\ot n} \ot W_1 ((v_1,\ldots ,v_n))[[h]]\quad\text{for all }n\geqslant 1.
$$
By substituting the variables $v=(v_1,\ldots ,v_n)$   with $ (z+u_1,\ldots z+u_n)$, for some variables $u=(u_1,\ldots ,u_n)$ and $z$,  we get, due to the expansion convention from Section \ref{sec13},
$$
Y_W(T_{[n]}^+ (u)\vac,z)w_1=\Tc_{[n]}(u|z)w_1\,\in\,\left(\ndo\CC^N\right)^{\ot n} \ot W_1 ((z))[[u_1,\ldots ,u_n,h]]
$$
for all $n\geqslant 1$.
Hence $W_1$ is a $\V_c(\R)$-submodule of $W$.
\end{prf}

\section*{Acknowledgement}
The author would like to thank Alexander Molev for useful discussions. 
We would also like to thank the anonymous referee for useful suggestions which helped us to improve the manuscript.
The research was partially supported by the  Croatian Science Foundation under the project 2634 and by the Australian Research Council.


\begin{thebibliography}{9}
\bibitem{C}
I. V. Cherednik, 
{\em A  new  interpretation  of  Gelfand--Tzetlin bases},
Duke Math. J. \textbf{54} (1987), 563--577.

\bibitem{D}
J. Ding, 
{\em Spinor Representations of $U_q(\hat{\gl} (n))$ and Quantum Boson-Fermion Correspondence},
Comm. Math. Phys. {\bf 200} (1999), 399--420;
\href{https://arxiv.org/abs/q-alg/9510014}{arXiv:q-alg/9510014}.

\bibitem{EK4}
P. Etingof and D. Kazhdan,
{\em Quantization of Lie bialgebras, IV},
Selecta Math. (N.S.) {\bf 6} (2000), 79--104;
\href{https://arxiv.org/abs/math/9801043}{arXiv:math/9801043 [math.QA]}.

\bibitem{EK}
P. Etingof, D. Kazhdan,
{\em Quantization of Lie bialgebras, V}, Selecta Math. (N.S.) \textbf{6} (2000), 105--130;
\href{http://arxiv.org/abs/math/9808121}{arXiv:math/9808121 [math.QA]}.

\bibitem{FJMR}
L. Frappat, N. Jing, A. Molev and E. Ragoucy,
{\em Higher Sugawara operators for the quantum affine
algebras of type $A$},
Comm. Math. Phys. {\bf 345} (2016), 631--657;
\href{https://arxiv.org/abs/1505.03667}{arXiv:1505.03667 [math.QA]}.

\bibitem{FZ}
I. B. Frenkel and Y.-C. Zhu,
{\em Vertex operator algebras associated to representations of affine and Virasoro algebras},
Duke Math. J. \textbf{66} (1992), 123--168.

\bibitem{GS}
D. Gurevich, P. Saponov,
{\em Centers in Generalized Reflection Equation algebras},
\href{https://arxiv.org/abs/1712.06154}{arXiv:1712.06154 [math.QA]}.

\bibitem{I}
K. Iohara,
{\em Bosonic representations of Yangian double $DY_{\hbar}(\mathfrak{g})$ with $\mathfrak{g}=\gl_N,\sll_N$},
J. Phys. A \textbf{29} (1996), 4593--4621;
\href{https://arxiv.org/abs/q-alg/9603033}{arXiv:q-alg/9603033}.

\bibitem{JKMY}
N. Jing, S. Ko\v{z}i\'{c}, A. Molev, F. Yang,
{\em Center of the quantum affine vertex algebra in type $A$},
 J. Algebra \textbf{496} (2018), 138--186;
\href{https://arxiv.org/abs/1603.00237}{arXiv:1603.00237 [math.QA]}.

\bibitem{JW}
D. Jordan, N. White,
{\em The center of the reflection equation algebra via quantum minors},
\href{https://arxiv.org/abs/1709.09149}{arXiv:1709.09149 [math.QA]}.

\bibitem{J}
A. Jucys,
{\em On the Young operators of the symmetric group},
Lietuvos Fizikos Rinkinys {\bf 6} (1966), 163--180.

\bibitem{Kas}
C. Kassel,
{\em Quantum Groups}, 
Graduate texts in mathematics; vol. \textbf{155}, Springer-Verlag, 1995.

\bibitem{Kho}
S. M. Khoroshkin,
{\em Central Extension of the Yangian Double},
\href{https://arxiv.org/abs/q-alg/9602031}{arXiv:q-alg/9602031}.

\bibitem{K}
S. Ko\v{z}i\'{c}, 
{\em Quasi modules for the quantum affine vertex algebra in type $A$}, 
Comm. Math. Phys. \textbf{365} (2019), 1049--1078;
\href{https://arxiv.org/abs/1707.09542}{arXiv:1707.09542 [math.QA]}.

\bibitem{KS}
P. P. Kulish and E. K. Sklyanin,
{\em Algebraic structures related to reflection equations}, 
J. Phys. A \textbf{25} (1992), 5963--5975;
\href{https://arxiv.org/abs/hep-th/9209054}{arXiv:hep-th/9209054}.

\bibitem{KJC}
V. B. Kuznetsov, M. F. J\o{rgensen}, P. L. Christiansen,
{\em New boundary conditions for integrable lattices},
J. Phys. A \textbf{28} (1995), 4639--4654;
\href{https://arxiv.org/abs/hep-th/9503168}{arXiv:hep-th/9503168}.

\bibitem{LLi}
J. Lepowsky, H.-S. Li,
{\em Introduction to Vertex Operator Algebras and Their Representations},
Progress in Math., Vol. \textbf{227}, Birkhauser, Boston, 2004.

\bibitem{Li_g1}
H.-S. Li, 
{\em Axiomatic $G_1$-vertex algebras}, 
Commun. Contemp. Math. \textbf{5} (2003), 281--327;
\href{https://arxiv.org/abs/math/0204308}{arXiv:math/0204308 [math.QA]}.

\bibitem{Li05}
H.-S. Li, 
{\em Nonlocal vertex algebras generated by formal vertex operators}, 
Selecta Math. (New Series) \textbf{11} (2005), 349--397;
\href{http://arxiv.org/abs/math/0502244}{arXiv:math/0502244 [math.QA]}.

\bibitem{Li06}
H.-S. Li, 
{\em Constructing quantum vertex algebras}, 
Int. J. Math. \textbf{17} (2006), 441--476;
\href{https://arxiv.org/abs/math/0505293v1}{ 	arXiv:math/0505293 [math.QA]}.

\bibitem{Li}
H.-S. Li,
{\em $\hbar$-adic quantum vertex algebras and their modules},
Comm. Math. Phys. {\bf 296} (2010), 475--523;
\href{http://arxiv.org/abs/0812.3156}{arXiv:0812.3156 [math.QA]}.

\bibitem{Lian}
B.-H. Lian,
{\em On the classification of simple vertex operator algebras},
Comm. Math. Phys. {\bf 163} (1994), 307--357.

\bibitem{MRS}
 M. Mintchev, E. Ragoucy  and P. Sorba,
{\em Spontaneous symmetry breaking in the $gl(N)-NLS$ hierarchy on the half line},
J. Phys. A \textbf{34} (2001) 8345--8364;
\href{https://arxiv.org/abs/hep-th/0104079}{arXiv:hep-th/0104079}.

\bibitem{M}
A. Molev,
{\em Yangians and classical Lie algebras}, 
Mathematical Surveys and Monographs, 143. American Mathematical Society, Providence, RI, 2007.

\bibitem{MR}
A. I. Molev, E. Ragoucy,
{\em Representations of reflection algebras},
Rev. Math. Phys. \textbf{14} (2002), 317--342;
\href{https://arxiv.org/abs/math/0107213}{arXiv:math/0107213 [math.QA]}.

\bibitem{O}
{A. Okounkov},
{\em Quantum immanants and higher Capelli identities},
{Transform. Groups} \textbf{1} (1996), 99--126;
\href{https://arxiv.org/abs/q-alg/9602028}{arXiv:q-alg/9602028}.

\bibitem{RS}
N. Yu. Reshetikhin, M. A. Semenov-Tian-Shansky,
{\em Central extensions of quantum current groups}, 
Lett. Math. Phys., {\bf 19} (1990), 133--142.

\bibitem{S}
E. K. Sklyanin, 
{\em Boundary conditions for integrable quantum systems}, 
J. Phys. A \textbf{21} (1988), 2375--2389.

\bibitem{T}
D. V. Talalaev,
{\em The quantum Gaudin system},
Funct. Anal. Appl. {\bf 40} (2006), 73--77.

\end{thebibliography}
\end{document}